\documentclass{amsart}

\usepackage{amssymb}

\usepackage{epic, eepic}
\def\squarebox#1{\hbox to #1{\hfill\vbox to #1{\vfill}}} 
\newcommand{\stopthm}{\hfill\hfill\vbox{\hrule\hbox{\vrule\squarebox 
                 {.667em}\vrule}\hrule}\smallskip} 

\newcommand{\CC}{{\mathbb C}}
\newcommand{\CI}{{\mathcal C}^\infty }
\newcommand{\CIc}{{\mathcal C}^\infty_{\rm{c}} }
\newcommand{\Op}{{\operatorname{Op}^{{w}}_h}}

\newcommand{\RR}{{\mathbb R}}

\newcommand{\ZZ}{{\mathbb Z}}
\newcommand{\MO}{{\mathcal M}}

\newcommand{\KKer}{\operatorname{ker}}
\newcommand{\supp}{\operatorname{supp}}

\newcommand{\UC}{U_{\circlearrowright}}
\newcommand{\QC}{Q_{\circlearrowright}}

\newcommand{\tr}{\operatorname{tr}}
\newcommand{\rest}{\!\!\restriction}

\renewcommand{\Re}{\mathop{\rm Re}\nolimits}
\renewcommand{\Im}{\mathop{\rm Im}\nolimits}

\theoremstyle{plain}

\newtheorem{thm}{Theorem}
\newtheorem{prop}{Proposition}[section]
\newtheorem{lem}[prop]{Lemma}

\theoremstyle{definition}

\numberwithin{equation}{section}

\def\bbbone{{\mathchoice {1\mskip-4mu \text{l}} {1\mskip-4mu \text{l}}
{ 1\mskip-4.5mu \text{l}} { 1\mskip-5mu \text{l}}}}

\title
{Quantum monodromy and semi-classical trace formul{\ae}} 
\author[J. Sj{\"o}strand]{Johannes Sj{\"o}strand}
\address{Centre de Math{\'e}matiques, {\'E}cole Polytechnique \\
UMR 7460, CNRS \\
F-91128 Palaiseau, France }
\email{johannes@math.polytechnique.fr}
\author[M. Zworski]{Maciej Zworski}
\address{Mathematics Department, University of California \\
Evans Hall, Berkeley, CA 94720, USA}
\email{zworski@math.berkeley.edu}

\begin{document}    
   
\maketitle   
   
\section{Introduction}   
\label{in}

Trace formul{\ae} provide one of the most elegant descriptions of the
classical-quantum correspondence. One side of a formula is given 
by a trace of a quantum object, typically derived from a quantum
Hamiltonian, and the other side is described in terms of closed 
orbits of the corresponding classical Hamiltonian. In algebraic
situations, such as the original Selberg trace formula, the identities
are exact, while in general they hold only in semi-classical or
high-energy limits. We refer to a recent survey \cite{Ur} for an
introduction and references.

In this paper we present an intermediate trace formula in which
the original trace is expressed in terms of traces of quantum monodromy
operators directly related to the classical dynamics. The usual trace
formul{\ae} follow and in addition this approach allows handling effective
Hamiltonians.

Let $ P = (1/i) h \partial_x $ be the semi-classical differentiation 
operator on the circle, $ x \in {\mathbb S}^1 = 
\RR / 2 \pi \ZZ $, $ 0 < h < 1 $. The 
classical Poisson formula can be written as follows: if $ \hat f \in 
{\mathcal C}^\infty _{\rm{c}} ( {\mathbb R} ) $ then 
\begin{equation}
\label{eq:poiss}
 {\mathrm{tr}} \; f ( P / h ) = \frac{1 }{ 2\pi i } 
\sum_{ |k|\leq N } \int_{\mathbb R} f (z/h) \left( e^{ 2 \pi i z / h } 
\right)^k \frac{d}{dz} \left( e^{ 2 \pi i z / h } \right) dz \,, 
\end{equation}
where $ N$ depends on the support of $ \hat f $, and we think of  
$M ( z , h ) = e^{ 2 \pi iz /h } : 
{\mathbb C} \rightarrow {\mathbb C} $ as the monodromy 
operator for the solutions of $ P - z $. It acts on functions in one dimension
lower (zero dimension here), identified geometrically with the functions 
on the transversal to the closed curve ($ {\mathbb S}^1 $ here), and 
analytically with $ {\mathrm{ker}}\; ( P - z ) $ (${\mathbb C} $ here).

Now let $ P $ be a semi-classical, self-andjoint, principal type operator,
with symbol $ p$ 
(for instance $ P = -h^2 \Delta + V ( x) $,  $ p = \xi^2 + V( x) $), and let
$ \gamma \subset p^{-1} (0 ) $ be a closed primitive orbit of the Hamilton flow of $ p$. 
We can define the {\em monodromy operator}, $ M ( z, h ) $ for $ P - z $ 
along $ \gamma$, acting on functions in one dimension lower, that is, 
on functions on the transversal 
to $ \gamma $ in the base. We then have 
\begin{thm}
\label{t:0}
Suppose that there exists a neighbourhood of $ \gamma $, 
$ \Omega $, satisfying the condition
\begin{equation}
\label{eq:cond}
 m \in \Omega  \ \text{and} \ \exp t H_p ( m ) = m \,, \ 
p ( m ) = 0 \,, \ \ 0 <  |t| \leq T N 
\ \Longrightarrow \ m \in \gamma \,, \end{equation}
where $ T $ is the primitive period of $ \gamma$. 
If $ \hat f \in 
{\mathcal C}^\infty_{\rm{c}} ( {\mathbb R} )$, 
$ {\rm supp}\; \hat f \subset ( - NT +  C , NT - C ) \setminus \{ 0\}$, $ 
C= C ( p ) \geq 0 $, 
$ \chi \in \CIc (\RR ) $, and $ A \in \Psi_h^{0,0} ( X ) $ is a microlocal 
cut-off to a sufficiently small neighbourhood of $ \gamma $, then
\begin{equation}
\label{eq:main}
 {\mathrm{tr}} \; f ( P / h ) 
\chi ( P ) 
A = \frac{1 }{ 2\pi i } 
\sum_{  -N - 1}^{N - 1 } {\mathrm {tr}} \; 
\int_{\mathbb R} f (z/h) M(z, h)^k \frac{d}{dz} M ( z , h )
 \chi ( z) 
dz + {\mathcal O} ( h^\infty ) \,, \end{equation}
where $ M ( z, h )$ is the semi-classical monodromy operator associated to 
$ \gamma $. 
\end{thm}
The dynamical assumption on the operator means
that in a neighbourhood of $ \gamma $ there are no other closed orbits of 
period less than $ T N$, on the energy surface $ p  = 0 $. We avoid 
a neighbourhood of $ 0 $ in the support of $ \hat f $ to avoid the
dependence on the microlocal cut-off $ A $.

The monodromy operator quantizes the Poincar{\'e} map for $ \gamma $ and 
its geometric analysis gives the now standard trace formul{\ae} of 
Selberg, Gutzwiller and Duistermaat-Guillemin (see \cite{CRR} for a 
recent proof and a historical discussion, and Sect.\ref{tr} for a derivation 
based on Theorem \ref{t:0}).
 The term $ k = -1 $ corresponds
to the contributions from ``not moving at all'' and the other terms to 
contributions from going $ |k+1| $ times around $ \gamma $, in the positive
direction when $ k \geq 0$, and in the negative direction, when $ k < -1 $.
For non-degenerate orbits we analyse the traces on monodromy operators in 
Sect.\ref{tr} and recover the usual semi-classical trace formul{\ae} in 
our general setting -- see Theorem \ref{t:2}.

Theorem \ref{t:0} is a special case of the more general Theorem \ref{t:1}
presented in Sect.\ref{pf}. Motivated by {\em effective Hamiltonians} in which 
the spectral parameter appears non-linearly, we give there a trace formula
for a family $ P ( z ) $ with the special case corresponding to 
$ P - z $. For an example of a use of effective Hamiltonians in an 
interesting physical situation we refer to \cite{dHvA}. The effective
Hamiltonian described there comes from the ``Peierls substitution'',
and the celebrated ``Onsager rule'' is a consequence of a calculation 
of traces.

The point of view taken here is purely semi-classical but when 
translated to 
the special case of $\CI$-singularities/high energy regime, it
is close to that of Marvizi-Melrose \cite{MM} and Popov \cite{Pop}. 
In those works the trace of the wave group was reduced to the study
of a trace of an operator quantizing the Poincar\'e map. In \cite{Pop}
it was used to determine contributions of degenerate orbits and our
formula could be used for that as well.

\medskip
\noindent
{\sc Acknowledgments.} We would like to thank the 
Erwin Schr\"odinger Institute, and the 
France-Berkeley Fund for 
providing partial support for this work. The second author
gratefully acknowledges the partial support of the National Science Foundation 
under the grant DMS-9970614.

\section{Outline of the proof}
\label{out}

To present the idea of our proof we use it to derive 
the classical Poisson
summation formula \eqref{eq:poiss}. The left hand side there
can be written using the usual functional calculus based on Cauchy's formula:
\begin{equation}
\label{eq:2.cau}
\tr f  \left( \frac{P}{ h } \right) =  \frac{1}{ 2 \pi i}
\tr \int_\Gamma f\left ( \frac{z}{h} \right) ( P - z )^{-1} dz \,, \ 
\ \Gamma = \Gamma_+ - \Gamma_- \,, \ \ \Gamma_{\pm } = \RR \pm i R \,,
\end{equation}
where we take the positive orientation of $ \RR $ and $ R > 0 $ is an 
arbitrary constant. We make an assumption on the support of the Fourier
transform on $ f$:
\begin{equation}
\label{eq:2.f}
\supp \hat f  \subset ( - 2 \pi N , 2 \pi N ) \,.
\end{equation}

We would like to replace $ ( P - z ) ^{-1} $
by an {\em effective Hamiltonian} which measures the obstruction
to the solvability of $ ( P - z ) u = f $. For that we introduce a
{\em Grushin problem} (see for instance \cite{Gr} for 
applications of this method in spectral problems, and for references):
\begin{equation}
\label{eq:2.gr}
{\mathcal P } ( z ) \stackrel{\rm{def}}{=} 
\left( \begin{array}{ll} P - z & R_- ( z) \\
R_+ ( z )  & 0 \end{array} \right) \; : \; H^1 ( {\mathbb S} ^1 ) \times \CC 
\ \longrightarrow \ L^2 ( {\mathbb S} ^1 ) \times \CC \,,
\end{equation}
where $ R_ \pm ( z)  $ should be chosen so that $ {\mathcal P} ( z) 
$ is invertible. 
If we put 
\[  R_+ u \stackrel{\rm{def}}{=} u ( 0 ) \,, \]
then we can locally solve
\[  \left\{ \begin{array}{l} ( P - z ) u = 0 \\
R_+ u = v \end{array} \right.  \,,  \] 
by putting 
\[ u = I_+( z)  v  = \exp ( i z x / h ) v \,, \ \ \  - \epsilon < x < 
2\pi - 2 \epsilon \,.\]
This is the forward solution, and we can also define the backward one
by 
\[ u = I_- ( z) v  = \exp ( i z x / h ) v \,, \ \ \  - 2 \pi  + 2 \epsilon < x 
 \epsilon \,.\]
The monodromy operator $ M ( z , h ) : \CC \rightarrow \CC $, can be 
defined by 
\begin{equation}
\label{eq:2.mon}  I_+ ( z) v ( \pi )  = I_- ( z ) M ( z , h )v  ( \pi ) 
 \,, \end{equation}
and we immediately see that
\[ M ( z , h ) = \exp \frac{ 2 \pi i  z }{ h } \,. \]
We use $ I_\pm ( z) $ and the point $ \pi $ to work with objects 
defined on $ {\mathbb S}^1 $ rather than on its cover: a more intuitive
definition of $ M ( z , h ) $ can be given by looking at a value of the
solution after going around the circle but that has some technical 
disadvantages.

Let $ \chi \in \CI ( {\mathbb S}^1 , [0 , 1 ] ) $  have the properties
\[ \chi ( x ) \equiv 1 \,, \ \ - \epsilon < x < \pi + \epsilon \,, \ \ 
\chi ( x) \equiv 0 \,, \ \   - \pi + 2 \epsilon < x < -2 \epsilon \,, \]
and put 
\[ E_+ ( z ) =  \chi I_+ ( z) + ( 1 - \chi ) I _- ( z ) \,.\]
We see that 
\[  ( P - z) E_+  =  [ P , \chi ] I_+ ( z) - [ P , \chi ] I_- ( z) 
=  [ P , \chi ]_- I_+ ( z) - [ P , \chi ]_- I_- ( z)  \,, \]
where $ [ P , \chi]_- $ denotes the part of the commutator supported near
$ \pi $.  This can be simplified using \eqref{eq:2.mon}:
\[  ( P - z ) E_+ + [ P , \chi] _- I_- ( z ) ( I - M ( z , h ) ) = 0 \,, \]
which suggests putting 
\[ R_- ( z ) = [ P , \chi]_- I_- ( z) \,, \]
so that the problem
\[ \left\{ \begin{array}{l} ( P - z) u + R_- ( z) u_ - = 0 \\
R_+ ( z ) u = v \end{array} \right. \]
has a solution:
\[ \left\{ \begin{array}{l} u = E_+ ( z ) v \\
u_- = E_{- + } ( z) v  \end{array} \right. \,, \]
with $ E_{-+} ( z) = I - M( z , h )\,. $ 

One can show\footnote{In this situation it is quite easy but it will
be done in greater generality in Sect.\ref{gr}.}
that with this choice of $ R_\pm ( z ) $, \eqref{eq:2.gr}
is invertible and then
\[ {\mathcal P}(z)^{-1} = {\mathcal E} ( z) = 
\left( \begin{array}{ll} E ( z ) & E_+ ( z) \\
E _- ( z ) & E_{-+ } ( z ) \end{array} \right) \,, \]
where all the entries are holomorphic in $ z $, and 
$ E_+ ( z) $, $ E_{ -+}( z) $, are as above. The operator $ E_{-+} ( z) $
is the effective Hamiltonian in the sense that its invertibility controls
the existence of the resolvent:
\begin{equation}
\label{eq:res}
 ( P - z )^{-1} = E ( z) - E_+ ( z) E_{ - + } ( z) ^{-1} E_- ( z) \,. 
\end{equation}
Inserting this in \eqref{eq:2.cau} and using the holomorphy of $ E ( z) $
gives
\[ \begin{split}
\tr f  \left( \frac{P}{ h } \right) & =   - \frac{1}{ 2 \pi i}
 \int_\Gamma f\left ( \frac{z}{h} \right) \tr  
 E_+ ( z) E_{ - + } ( z) ^{-1} E_- ( z)  dz \\ & =
-   \frac{1}{ 2 \pi i } \int_\Gamma f\left ( \frac{z}{h} \right) \tr  
E_- ( z) E_+ ( z) E_{ - + } ( z) ^{-1}  dz  \,, \end{split} \] 
where we used the cyclicity of the trace.
Differentiating $ {\mathcal E} ( z ) {\mathcal P } ( z) = Id $ shows
\[ E_- ( z) E_+ ( z) = \partial_z E_{-+}(z) + E_- ( z ) \partial _z
R _- ( z) E_{-+} ( z ) \,, \]
which inserted in the previous identity gives
\[ \tr f  \left( \frac{P}{ h } \right) =  - \frac{1}{ 2 \pi i}
 \int_\Gamma f\left ( \frac{z}{h} \right) \tr  
\partial_z E_{- +} ( z)  E_{ - + } ( z) ^{-1}  dz \,, \]
where we eliminated the other term using countour deformation.

We now use the expression for $ E_{ - + } $ to write 
\[ \begin{split} \tr f  \left( \frac{P}{ h } \right) = &   \frac{1}{ 2 \pi i}
 \int_{\Gamma_+}  f\left ( \frac{z}{h} \right) \tr  
\partial_z M( z , h )   ( I - M ( z,h ) )  ^{-1}  dz  \\
&  + \;  
\frac{1}{ 2 \pi i } 
\int_{\Gamma_{ -} } 
  f\left ( \frac{z}{h} \right) \tr  
\partial_z M( z , h ) M ( z, h )^{-1}  ( I - M ( z,h )^{-1} )  ^{-1}  dz 
\,.   \end{split} \]

The assumption \eqref{eq:2.f} and the Paley-Wiener theorem give

\[ | \hat f ( z / h ) | \leq e^{ 2 \pi N | \Im z |/h  } \langle \Re z  / h 
\rangle^{ -\infty } \,. \]
Writing 
\[   ( I - M ( z,h ) )  ^{-1} = \sum_{ k = 0 } ^{ N -1} M ( z , h ) ^ k 
+ M ( z , h ) ^{ N} ( I - M ( z , h ))^{-1} \,, \]
for $ \Gamma_+ $, and 
\[  M ( z, h )^{-1}  ( I - M ( z,h )^{-1} )  ^{-1} =
 \sum_{ k = 1 } ^{ N} M ( z , h ) ^ {-k }
+ M ( z , h ) ^{ - N - 1} ( I - M ( z , h ))^{-1} \,, \]
for $ \Gamma_- $, we can eliminate the 
last terms by deforming the contours to imaginary infinities ($  R 
\rightarrow \infty $), and this gives \eqref{eq:poiss}.

In the general situation we proceed similarly but now {\em microlocally}
in a neighbourhood of the closed orbit described in Theorem \ref{t:0}
-- see Sect.\ref{sc} for a precise definition of microlocalization.
The formula \eqref{eq:2.cau} has to be replaced by 
\begin{equation}
\label{eq:2.cau'}
\tr f  \left( \frac{P}{ h } \right) \chi ( P ) A = -  \frac{1}{ \pi }
 \int_{\CC}  f\left ( \frac{z}{h} \right) \bar \partial_z
\tilde \chi ( z) ( P - z )^{-1} A\;  {\mathcal L}( dz)  \,, 
\end{equation}
where $ \tilde \chi $ is an {\em almost analytic extension} of
$ \chi $, that is an extension satisfying $ \bar \partial_z \chi ( z) 
= {\mathcal O} ( | \Im z |^{\infty } ) $ -- see Sect.\ref{sc},
and we want to procceed with a similar reduction to the effective 
Hamiltonian given in terms of the {\em monodromy operator}.

To construct the monodromy operator we fix two different points on $ \gamma$,
$ m _0 $, $ m_1$ (corresponding to $ 0 $ and $ \pi $ in the example),
and their disjoint neighbourhoods, $ W_+ $ and $ W_- $ respectively.
We then consider local kernels of $ P - z$ near $ m_0 $ and $ m_1 $ (that
is, sets of disctributions satisfying $ ( P -z ) u = 0 $ near $ m_i $'s),
$ \KKer_{m_j} ( P - z ) $, $ j = 0, 1$,  with elements microlocally 
defined in $ W_\pm $.
and the forward and backward solutions:
\[ I_\pm ( z) \; : \; \KKer_{ m_0 } ( P - z) \ \longrightarrow \
\KKer_{ m_1 } ( P - z ) \,.  \]
We then define the {\em quantum monodromy operator}, $ {\mathcal M} (z)$
by
\begin{gather*}
 I _{-}  ( z ) {\mathcal M} ( z ) = I_+ ( z ) \,, \ \ \ \
{\mathcal M} ( z) \; : \;   \KKer_{ m_0 } ( P - z)
\ \longrightarrow \  \KKer_{ m_0 } ( P - z) \,. \end{gather*}
The operator $ P $ is assumed to be self-adjoint with respect to some
inner product $ \langle \bullet, \bullet \rangle$, and we define
the {\em quantum flux} norm on $ \KKer_{m_0} ( P - z) $ as 
follows\footnote{See \cite{Gr} for an earlier mathematical 
development of this basic quantum mechanical idea.}: let $ \chi $ be 
a microlocal cut-off function, with basic properties of the function
$ \chi $ in the example. Roughly speaking $ \chi $ should supported near
$ \gamma $ and be equal to one near the part of $\gamma $ between $ W_+$ and 
$ W_- $. We denote by $ [ P , \chi ]_{W_+} $ the part of the
commutator supported in $ W_+$, and  put
\[  \langle u , v \rangle_{\rm{QF}}  \stackrel{\rm{def}}{=} 
\langle [ ( h / i ) P , \chi ]_{W_+} u , v \rangle \,, \ \ u ,v 
\in \KKer_{m_0} ( P - z ) \,. \]
It is easy to check that this norm is independent of the choice of 
$ \chi $ -- see the proof of Lemma \ref{l:4.1}.  This independence leads
to the unitarity of $ {\mathcal M} ( z ) $:
\[  \langle {\mathcal M} ( z)  u ,  {\mathcal M} ( z )  u \rangle_{\rm{QF}}
=  \langle u , u  \rangle_{\rm{QF}} \,, \ \ u \in \KKer_{m_0} ( P - z )  
\,. \]
For practical reasons we identify $ \KKer_{ m_0 } ( P - z) $ with 
$ {\mathcal D}' ( \RR^{n-1} ) $, microlocally near $ ( 0 , 0 )$, and
choose the idenfification so that the corresponding monodromy map is
unitary (microlocally near $ ( 0 , 0 ) $ where $ ( 0, 0 )$ corresponds
to the closed orbit intersecting a transversal identified with 
$ T^* \RR^{n-1} $). This gives
\[ M ( z , h ) \; : \; {\mathcal D} ' ( \RR^{n-1} ) \ \longrightarrow
\ {\mathcal D}' ( \RR^{n-1} ) \,, \]
microlocally defined near $ ( 0, 0 ) $ (see Sect.\ref{sc} for a precise
definition of this notion) and unitary there.
This is the operator appearing in Theorem \ref{t:0} and it shares
many properties with its simple version $ \exp(2 \pi i z / h ) $ appearing
for $ {\mathbb S}^1$. 

As shown in the example of the Poisson formula, traces can be expressed
in terms of traces of effective Hamiltonians ($ E_{-+} ( z) $ there).
Hence in our final formula, 
we replace $ P - z $ by a more general operator $ P ( z ) $, 
for which we do not demand holomorphy $ z $ but
only that $ P ( z) $ is self-adjoint for $ z $ real and that it
is an almost analytic family of operators. In Theorem \ref{t:1}
in Sect.\ref{pf} we will compute the trace of 
\[  - \frac1\pi  {\mathrm{tr}} \;  \int f ( z / h ) 
\;  \bar \partial_z  \left[ \tilde \chi ( z )\;  
 \partial_z P (z) \;  P ( z) ^{-1} \right]  A  \; 
{\mathcal L} ( d z )  \,, \]
which for $ P ( z) = P - z$ reduces to \eqref{eq:2.cau'}.

The only prerequisite to reading the paper is the basic calculus of
semi-classical pseudodifferential operators (see \cite{DiSj}).
In Sect.\ref{sc} we review
various aspects of semi-classical microlocal analysis needed here.
In Sect.\ref{qt} we define the {\em quantum time} and {\em quantum 
monodromy}. Then in Sect.\ref{gr} we follow the procedure described
for $ {\mathbb S}^1 $  to solve a Grushin problem allowing us to 
represent $ P ( z ) ^{-1} $ near a closed orbit. That is applied in the
proof of the trace formula in Sect.\ref{pf}, and in Sect.\ref{tr}
we derive the more standard trace formula in the case of a non-degenerate
orbit.

\section{Semi-classical operators and their almost analytic extensions}
\label{sc}

Let $ X $ be a compact $\CI$ manifold. 
We introduce the usual class of semi-classical 
symbols on $ X $:
\[ S^{m,k} ( T^* X ) = \{ a \in \CI( T^* X \times (0, 1]  ) :
|\partial_x ^{ \alpha } \partial _\xi^\beta a ( x, \xi ;h ) | \leq
C_{ \alpha, \beta} h^{ -m } \langle \xi \rangle^{ k - |\beta| } 
\} \,, \]
and the class  corresponding pseudodifferential operators, $ 
\Psi_{h}^{m,k} ( X ) $, with the quantization and symbol maps:
\[ 
\begin{split}
& \Op \; : \; S^{m , k } ( T^* X ) \ \longrightarrow 
\Psi^{m,k}_h ( X) \\
& \sigma_h \; : \; \Psi_h^{m,k} ( X ) \ \longrightarrow 
S^{ m , k } ( T^* X ) / S^{ m-1, k-1} ( T^* X ) \,, \end{split}
\]
with both maps  surjective,  and the usual properties 
\[ \sigma_h ( A \circ B ) = \sigma_h ( A)\sigma _h ( B ) \,,\]
\[ 0 \rightarrow \Psi^{m-1, k-1} ( X) \hookrightarrow \Psi^{m, k} ( X)
\stackrel{\sigma_h}{\rightarrow} S^{ m , k } ( T^* X ) 
/ S^{ m-1, k-1} ( T^* X ) \rightarrow 0 \,, \]
a short exact sequence, and
\[ \sigma_h \circ \Op : S^{m,k} ( T^* X ) 
\ \longrightarrow  S^{ m , k } ( T^* X ) / S^{ m-1, k-1} ( T^* X ) \,,
\]
the natural projection map.
The class of operators and the quantization map are defined locally using the
definition on $ \RR^n $:
\begin{equation}
\label{eq:weyl}
 \Op (a) u ( x) = \frac{1}{ ( 2 \pi h )^n } 
\int \int  a \left( \frac{x + y }{2}  , \xi \right
) e^{ i \langle x -  y, \xi \rangle / h } u ( y ) 
dy d \xi \,, \end{equation}
and we refer to \cite{DiSj} or \cite{WZ} for a detailed discussion.
We remark only that unlike the invariantly defined 
symbol map, $ \sigma_h $, 
the quantization map $ \Op$ can be chosen in many different ways.

In this paper we consider pseudo-differential operators as acting
on half-densities and consequently the symbols will also be considered
as half-densities -- see \cite[Sect.18.1]{Hor} for a general introduction
and the Appendix to this paper for a semi-classical discussion. 
For notational simplicity we supress the half-density notation. 
The only result we will need here is that in Weyl quantization, the
symbol is well defined up to terms of order $ {\mathcal O} (h^2 )$ 
-- see Appendix.


\medskip

For $ a \in S^{m,k} ( T^* X ) $ we define 
\[ {\text{ess-supp}}_h\; a \subset T^* X \sqcup S^* X \,, \ \ 
S^* X \stackrel{{\rm{def}}}{=}  ( T^* X  \setminus 0 ) / {\mathbb R}_+ \,, 
\]
where the usual $ \RR_+ $ action is given by multiplication on 
the fibers: $ ( x, \xi ) \mapsto  ( x , t \xi ) $, as 
\[
\begin{split}
&  {\text{ess-supp}}_h\; a =  \\
& \ \ \ \complement \{ ( x, \xi ) \in T^* X \; : \; 
\exists \; \epsilon > 0 \ \partial_x ^\alpha \partial_\xi ^\beta 
a ( x', \xi' ) = {\mathcal O} ( h^\infty ) 
  \,, \  d( x, x' ) + | \xi - \xi' |  < \epsilon \} \\
& \ \ \ \cup \complement  \{ 
( x, \xi ) \in T^* X \setminus 0 \; : \; 
\exists \; \epsilon > 0 \ \partial_x ^\alpha \partial_\xi ^\beta 
a ( x', \xi' ) = {\mathcal O} ( h^\infty  \langle \xi' \rangle^{  -\infty}) 
  \,,  \\ 
& \ \ \ \ \ \ \ \ \ 
d( x, x' ) + 1 / |\xi'| + | \xi/ |\xi|  - \xi'/ |\xi'| |  
< \epsilon \} / {\mathbb R}_+ 
\end{split} \,. \]
For $ A \in \Psi^{m,k} _h ( X) $, then define
\[ WF_h ( A) =  {\text{ess-supp}}_h\; a \,, \ \ A = \Op ( a ) \,, 
\]
noting that, as usual, the definition does not depend on the choice of
$ \Op $. For 
\[  u \in \CI ( ( 0 , 1]_h ; {\mathcal D}' ( X) ) \,, \ \  \exists \; 
N_0\,,  \ \ h^{- N_0} u \ \text{ is bounded in $ {\mathcal D}' ( X ) $,}\]
 we 
define
\[ WF_h ( u ) = \complement \{ ( x, \xi ) \; : \; 
\exists \;  A \in \Psi_h^{0,0} ( X) \ \sigma_h ( A ) ( x, \xi ) \neq 0 \,,
\ A u \in h^\infty  \CI ( ( 0 , 1]_h ; \CI ( X) ) \} \,.\]
When $ u $ is not necessarily smooth we can give a definition 
analogous to that of $ {\text{ess-supp}}_h\; a $.
Since in this note we will be concerned with a purely semi-classical 
theory and deal only with {\em compact} subsets of $ T^* X $ this 
definition is sufficient for our purposes (for more general definitions
of wave front set which include this usual semi-classical definition,
see \cite{MZ}).

To discuss {\em almost analytic continuation} of semi-classical 
pseudodifferential operators let us first recall the scalar case.
For $ f \in \CI ( \RR ) $, an almost analytic extension of $ f $ is
$ \tilde f \in \CI ( \CC ) $ such that locally uniformly
\[ \bar \partial_z \tilde f ( z ) = {\mathcal O} ( |\Im z |^\infty ) 
\,, \ \ \tilde f \rest_{\RR} = f \,. \]
The almost analytic extensions were introduced by H{\"o}rmander 
and are unique up to 
$ {\mathcal O} ( |\Im z|^\infty ) $ terms (see \cite[Sect.8]{DiSj}
and references given there).

Suppose now that 
\[ A ( x ) \in \CI ( \RR_x ; \Psi^{m,k} _h ( X) ) \]
is a smooth family of pseudodifferential operators. We can then find
$ a ( x ) \in \CI ( \RR_x ; S^{ m , k } ( T ^* X ))  $ such that $ 
A ( x ) = \Op ( a ( x) ) $. We then define the {\em almost analytic
extension of the family} $ A ( x ) $ as
\[ \widetilde A ( z ) = \Op ( \tilde a ( z ) ) \,, \]
where $ \tilde a ( z ) \in \CI ( \CC_z ;  S^{ m , k } ( T ^* X ) ) $
is an almost analytic extension of $ a ( x )$. To justify this 
definition we need the following easy
\begin{lem}
\label{l:2.1}
If $ a ( x ) \in \CI ( \RR_x ; S^{ m , k } ( T ^* X ))  $ then there
exists an almost analytic extension of $ a $ satisfying
\[  \tilde a ( z ) \in \CI ( \CC_z ;  S^{ m , k } ( T ^* X ) ) \,, \ \ 
\partial_x^\alpha \partial_\xi^\beta \bar \partial
_z \tilde a ( z )( x, \xi ; h )
= {\mathcal O} ( | \Im z |^\infty \langle \xi \rangle^{ k - |\beta| } )\,.
\]
\end{lem}

We will also need certain aspects of the theory of semi-classical
Fourier Integral Operators. Rather than review the full theory we
will consider a special class, to which the general calculus reduces
in local situations. Thus let $ A ( t ) $ be a smooth 
family of pseudodifferential
operators, $ A ( t ) = \Op ( a ( t ) ) $, 
$ a (t)  \in \CI ([-1  , 1]_t ; S^{ 0, -\infty }
( T^* X ))  $, such
that for all $ t $, $ WF ( A ( t ) ) \Subset T^* X $. We then define
a family of operators
\begin{gather}
\label{eq:2.2}
\begin{gathered}
  U( t ) \; : \; L^2 ( X) \rightarrow L^2 ( X) \,, \\
h D_t U ( t) +  U ( t) A(t) = 0 \,, \ \ U( 0 ) = U_0 \in \Psi_h^{0,0} ( X) \,.
\end{gathered}
\end{gather}
This is an example of a family of $h$-Fourier Integral Operators, $ U ( t)$,
associated to canonical transformations $ \kappa(t)$, 
generated by the Hamilton 
vector fields $ H_{a_0(t)} $, where the real valued
$ a_0 (t) $ is the $h$-principal 
symbol of $ A ( t ) $,
\[ \frac{d}{dt} \kappa( t) ( x, \xi ) =  ( \kappa ( t ))_* ( H_{a_0 ( t)} 
( x, \xi )) \,, \ \ \kappa( 0 ) ( x, \xi ) = ( x, \xi )\,,   \ \ 
( x, \xi ) \in T^* X \,.\]
All that we will need in this note is the {\em Egorov theorem} which 
can be proved directly from this definition: 
when $ U_0 $ in \eqref{eq:2.2} is elliptic (that is 
$ | \sigma (U_0 )| > c > 0 $  on $ T^* X $, then 
for $ B \in \Psi_h^{m,k} ( X ) $
\begin{gather}
\label{eq:egorov}
\begin{gathered}
\sigma ( V(t) B U ( t) ) = ( \kappa ( t) )^* \sigma ( B ) \,, \\
 V(t) U ( t ) - I \,, \,U ( t) V ( t) - I 
\in \Psi_h^{- \infty , - \infty } ( T^* X) 
\end{gathered}
\end{gather}
where the approximate inverse is constructed by taking
\[ h D_t V ( t) - A ( t) V(t) = 0 \,,  \ \ V ( 0 ) = V_0 \,, \ \ 
V_0 U_0 - I \,, \, U_0 V_0 - I \in \Psi_h^{- \infty , - \infty } ( T^* X) 
\,, \]
the existence of $ V_0 $ being guaranteed by the ellipticity of $ U_0 $.
The proof of \eqref{eq:egorov} follows from writing $ B(t) = 
V ( t) B U ( t) $, so that, in view of the properties of $ V (t)$, 
\[ h D_t B( t) \equiv [ A( t) , B( t ) ] \, \mod \Psi_h^{-\infty, -\infty}
\,, \ \ B( 0 ) = B_0 \,. \]
Since the symbol of the commutator is given by $ (h/i) H_{ a_0 ( t) }
\sigma ( B ( t) )$, \eqref{eq:egorov} follows directly from the definition
of $ \kappa( t)$.

If $ U = U(1) $, say, and the graph of $ \kappa ( 1) $ is denoted by $C$,
we conform to the usual notation and write
\[  U \in I^0_h ( X \times X ; C' ) \,, \ \ 
C' = \{ ( x, \xi;  y , - \eta) \; : \; ( x, \xi) = \kappa ( y , \eta )   \}
\,, \]
which means that $ U $ is an $h$-Fourier Integral Operator associated to 
the canonical graphs $ C$. Locally all $h$-Fourier Integral Operators
associated to canonical graphs are of the form $ U ( 1) $ thanks to the 
following well known 
\begin{lem}
\label{l:well}
Suppose that $ U_1, U_2  $ are open neighbourhoods of $ ( 0 , 0 ) 
\in T^* \RR^n $, and $ \kappa \; : \; U_1 \rightarrow U_2 $ is a 
canonical transformation satisfying $ \kappa ( 0 , 0 ) = ( 0 , 0 ) $.
Then there exists a smooth family of canonical transformations 
$ \kappa_t \; : \; U_1 \rightarrow U_2 $, $ 0 \leq t \leq 1 $,
satisfying $ \kappa_0 = id $, $ \kappa_1 = \kappa $, $\kappa_t ( 0, 0 ) 
= ( 0, 0 ) $. 
\end{lem}
\begin{proof}
Since the symplectic group, $ Sp ( n , \RR )$, is connected we can
first deform $ \kappa $ so that $ d \kappa ( 0 , 0 ) = Id$.  Hence,
near $ ( 0 , 0 ) $, $( ( x ( \kappa ( y , \eta )) , \xi ( \kappa ( y  ,
\eta ) ) ; y , \eta ) \mapsto ( x , \eta ) $ is surjective, and on the
graph of $ \kappa $, $ y $ and $ \eta $ can be regarded as functions
of $ x $ and $ \eta $. Since the symplectic forms, $ - d ( \langle y ,
d \eta \rangle )$ $ d ( \langle \xi , d x \rangle ) $ are equal, their
difference can be written locally as a differential:
\[   \langle y , d \eta \rangle  +     \langle \xi , d x \rangle 
= d \phi \,, \ \ \phi =  \phi ( x , \eta ) \,, \ \ d \phi ( 0 , 0 ) = 0 \,,\]
so that $ \kappa \; : \; ( \phi'_\eta ( x , \eta ) , \eta ) \mapsto ( 
x , \phi_x ' ( x , \eta ) ) $.  We could now take as our family
\[ \kappa_t \; : \; ( t \phi'_\eta ( x , \eta ) + ( 1 - t ) x , \eta ) \mapsto ( 
x , t \phi_x ' ( x , \eta ) +  ( 1 - t ) \eta )  \,. \]
The two steps can be connected smoothly by making the deformations 
flat at their junction. 
\end{proof}

The almost analytic continuation of a family of $h$-Fourier Integral 
Operators defined by \eqref{eq:2.2} is obtained by means of the following
\begin{lem}
\label{l:2.2}
Suppose that $ U ( t )$ is defined by \eqref{eq:2.2} and that $ 
\widetilde A ( z ) $ is an almost analytic continuation of the family 
$ A ( t )$, as given by Lemma \ref{l:2.1}.  Let $ \widetilde U ( 
z ) = \widetilde U ( t + i s ) $ be the solution of 
\begin{equation}
\label{eq:2.3}
\frac{1}{i} h D_s \widetilde U ( t + i s ) +  
 \widetilde U ( t + i s) \widetilde A ( t + i s ) = 0 \,, \ \ 
\widetilde U  ( t+ i s ) \rest_{ s = 0 } = U ( t ) \,.
\end{equation}
Then for $ |\Im z | \leq h \log h^{-L} $ we have 
\begin{eqnarray}
& & \| \widetilde U ( z ) \|_{L^2 \rightarrow L^2 } \leq  C 
\exp ( C | \Im z| /  h ) 
\label{eq:2.4}\\
& & \| \bar \partial_z \widetilde U ( z) \|_{ L^2 \rightarrow L^2 } 
= {\mathcal O} ( |\Im z |^\infty ) 
\label{eq:2.5}\\
& & h D_z \widetilde U ( z ) = \widetilde A (  z) \widetilde U ( z ) 
+ {\mathcal O}_{ L^2 \rightarrow L^2 } ( |\Im z|^\infty ) \,.
\label{eq:2.6}
\end{eqnarray}
\end{lem}
\begin{proof}
To see \eqref{eq:2.4} we write 
\[ h \frac{d}{ds} \| \widetilde U ( t + i s ) v \|^2 = 
2 \Re \langle   \widetilde U ( t + i s ) \widetilde A ( t + is ) v , 
\widetilde U ( t + i s)v \rangle )\leq C \| \widetilde U ( t + i s ) 
\|^2 \| v \|^2 \,, \ \ s > 0 \,. \]
Let us now take $ v $ with $ \| v \|=1 $ so that, by integration, 
\[ \| \widetilde U ( t + i s ) v \|^2  \leq \| \widetilde U ( t ) \|^2 
+ \frac{C}{h} \int_0^s \| \widetilde U ( t+ i \sigma ) \|^2 d \sigma \,.
\]
Since this holds for every $ v $ with $ \| v \| = 1 $ we can replace
the left hand side of the inequality by $ \| \widetilde U ( t + i s) \|^2 $,
and the standard Gronwall inequality argument shows that
\[ \| \widetilde U ( t + i s )  \|^2   \leq C e^{ Cs / h } 
 \,, \]
which is the desired bound.
Putting $ \widetilde V ( t + i s) = \bar \partial_z \widetilde U ( t + 
i s) $ we have 
\begin{gather*}
 h \partial_s \widetilde V ( t + i s ) = 
\widetilde U  ( t + i s ) 
 \bar  \partial_z \widetilde A  ( t + i s )  + 
\widetilde V  ( t + i s)  \widetilde A  ( t + i s ) =  
\widetilde V  ( t + is ) \widetilde  A 
 ( t + i s ) + {\mathcal O} ( s^\infty ) \,, \\ 
|s| < h \log h^{-L} \,, \  \widetilde V ( t + i s) \rest_{ s = 0 } = 0 \,,
\end{gather*}
where the initial condition came from the equation on the real axis:
$ h D_t \widetilde U ( t + i s ) \rest_{ s = 0 } =  U ( t )A ( t) 
 $. As in the argument for \eqref{eq:2.4},
this implies \eqref{eq:2.5} and \eqref{eq:2.6}.
\end{proof}

Our definitions of pseudo-differential operators and of (the special class of)
$h$-Fourier Integral Operators were global. It is useful and natural to 
consider the operators and their properties microlocally.
We consider classes of {\em tempered} operators:
\[ T  \; : \; \CI ( X ) \rightarrow \CI ( X) \,, \]
and for any semi-norms $ \| \bullet\|_1 $ and
$ \| \bullet\|_2 $ on $ \CI ( X) $ there exists $ M_0$ such that
\[  \| T u \|_1  = {\mathcal O } ( h ^{ - M_0 } ) 
\| u \|_2 \,. \]
For open sets, $ V \subset T^* X$, $ U \subset T^* X $, the operators 
{\em defined microlocally} near  $ V \times U $ 
are given by equivalence classes of tempered operators
given by the relation
\[ T \sim T' \ \Longleftrightarrow \ A ( T - T' ) B = {\mathcal O} ( h^\infty )
\; : \; {\mathcal D}' ( X) \ \longrightarrow \CI ( X ) \,, \]
for any 
$ A, B \in \Psi_h ^{0,0 } ( X ) $ such that 
\begin{gather}
\label{eq:2.7}
\begin{gathered}
WF ( A ) \subset \widetilde V \,, \ \ WF ( B ) \subset \widetilde U\,,
\\
\bar  V \Subset  \widetilde V \Subset T^* X \,, \ \ 
\bar U \Subset  \widetilde U \Subset T^* X \,, \ \ \widetilde U \,,\,
\widetilde V \ \text{ \  open}\,.
\end{gathered}
\end{gather}
The equivalence class $ T $, 
 $h$-Fourier Integral Operator associated to a local canonical graph
$ C $ if, again for any $ A $ and $ B$ above 
\[ A T B 
\in I^0 ( X \times X ; \widetilde C' ) \,,\]
where $ C$ needs to be defined only near $ U \times V $. 

We say that $ P = Q $ {\em microlocally} near $ U \times V $ if 
$  A P B - A Q B = {\mathcal O}_{L^2 \rightarrow L^2 }
 ( h ^ \infty ) $, where because of the assumed pre-compactness of 
$ U $ and $ V$ the $ L^2 $ norms can be replaced by any other norms.
For operator identities this will 
be the meaning of equality of operators in this paper, with $ U , V $
specified (or clear from the context). Similarly, we say that $ 
B = T^{-1} $ microlocally near $ V \times V $, if $ B T = I$ microlocally 
near $ U \times U $, and $ T B = I $ microlocally near $ V \times U $.
More generally, we could say that $ P = Q $ microlocally on $ 
W \subset T^* X \times T^* X $ (or, say, $ P $ is 
microlocally defined there), if for any $ U, V$, $ U \times V \subset
W $, $ P = Q $ microlocally in $ U \times V$. 
We should stress that ``microlocally'' is always meant in this 
semi-classical sense in our paper.

In this terminology we have a 
characterization of local $h$-Fourier Integral operators,
which is essentially the converse of Egorov's theorem:
\begin{lem}
\label{l:egor}
Suppose that $ U = {\mathcal O} ( 1) : L^2 ( X) \rightarrow 
L^2 ( X ) $, and that for every $ A \in \Psi_h ^{0,0} ( X ) $ we have
\[ A U = U B \,, \ \ B \in \Psi^{0,0}_h ( X) \,, \ \ 
\sigma ( B ) = {\kappa}^* \sigma ( A ) \,,\]
microlocally near $ ( m_0 , m_0 )$ 
where $ \kappa : T^* X \rightarrow T^*X $ is a symplectomorphism,
defined locally near $ m_0 $, $ \kappa ( m_0 ) = m_0 $. 
Then 
\begin{gather*}
 U \in  I^0_h ( X \times X ; C' ) \,, \ \text{ microlocally near 
$ ( m_0 , m_0 ) $, } \\ 
C' = \{ ( x, \xi;  y , - \eta) \; : \; ( y , \eta ) = \kappa ( x, \xi) \}
\,.\end{gather*}
\end{lem}
\begin{proof}
From Lemma \ref{l:well} we know that there exists a family of local 
symplectomorphisms, $ \kappa_t $, satisfying $ \kappa_t ( m_0 ) = m_0 $,
and $ \kappa_1 = \kappa $, $ \kappa_0 = id $. Since we are working
locally, there exists a function $ a(t) $, such that $ \kappa_t $ is
generated by its Hamilton vectorfield $ H_{a(t)} $. Let us now consider
\[ hD_t U(t) =  U ( t) A ( t) \,, \ \ U( 1) = U \, , \ \ 0 \leq t \leq 1
\,. \]
The same arguments as the one used in the proof of \eqref{eq:egorov}
shows that $ U(0)$ satisfies 
\begin{equation}
\label{eq:comm}
 [ U ( 0 ) , A ] = {\mathcal O} ( h ) \,, \ \ \text{
for any $ A \in \Psi^{0,0}_h ( X ) $. } \end{equation}
In fact, we take $ V (t ) $ with $ V( 0 ) = Id $ microlocally near $
( m_0 , m_0 ) $, so that 
\[  A U ( t ) V ( t ) = U ( t ) V( t) ( 
V(t)^{-1} B V ( t) ) = U ( t) V( t) A + {\mathcal O} ( h ) \,,
\]
where we
used Egorov's theorem and the assumption that $ \sigma ( B ) = 
\kappa^* \sigma ( A) $.  Putting $ t = 0 $ gives \eqref{eq:comm}.
By Beals's Lemma 
\cite[Prop.8.3]{DiSj} we conclude that $ U ( 0 ) \in \Psi^{0,0}_h ( X ) $,
and hence $ U $ is a microlocally defined $h$-Fourier Integral Operator 
associated to $ \kappa $.
\end{proof}

If the open sets 
$ U $ or $ V $ in \eqref{eq:2.7}
are small enough, so that they can be identified with 
neighbourhoods of points in $ T^* \RR^n $, we can use that identification to 
state that $ T $ is microlocally defined near, say, $ ( m , ( 0 , 0 ) )$, 
$ m \in T^* X $, $ ( 0, 0 ) \in T^* \RR^n $. An example useful here is
given in the next proposition.

By Darboux's theorem we know that
if $ p $ is a function with a non-vanishing differential then
there exists a local canonical transfomation $ \kappa $ such that 
$ \kappa^* p = \xi_1 $ where $ \xi_1 $ is part of a coordinate system
in which the symplectic form is the canonical one $ d ( \langle \xi , dx
\rangle ) $. The quantum version is given in
\begin{prop}
\label{p:2.1}
Suppose that $ P \in \Psi_{h}^{0,k} ( X ) $ is a semi-classical 
{\em real principal type operator}: $ p = \sigma ( P ) $ is real, 
independent of $ h$, and  
\[ p = 0 \Longrightarrow dp \neq 0   \,.\]
For any $ m_0 \in \Sigma_p 
\stackrel{\rm{def}}{=}
\{ m \in T^* X \; : \; p ( m ) = 0 \} $
 there exists a local canonical transformation
$ \kappa : T^* X \rightarrow T^* \RR^n $ defined near $ (  ( 0, 0 ), m_0 ) $, 
and an $h$-Fourier Integral Operator,  $ T $, 
associated to its graph, such that
\begin{gather*}
\kappa^*  \xi_1 = p \,,\\
T P = hD_{x_1}  T \,, \ \ \text{microlocally near $( ( 0 , 0 ),  m_0  )$}\, \\
T^{-1} \ \ \text{exists microlocally near $(m_0 , ( 0 , 0 ) ) $}\,.
\end{gather*}
\end{prop}
For the reader's convenience we outline a self-contained proof of this
semi-classical analogue of the standard $\CI$ result 
\cite[Proposition 26.1.3$'$]{Hor}.
\begin{proof}
By assumption $ dp ( m_0 ) \neq 0 $, and consequently Darboux's theorem gives
$ \kappa $ with the desired properties. Lemma \ref{l:well} then gives
us a family of symplectic transfomations $ \kappa_t $. 
If $ T_0 = U ( 1) $, where $ U ( 1) $ was defined using the family 
$ \kappa_t $, then \eqref{eq:egorov} shows that
$ T_0 P - hD_{x_1} = E  \in \Psi^{-1, 0} $ microlocally near $ ( 0, 0 )$. 
Hence we look for $ A$ such that 
$  hD_{x_1} + E = A h D_{x_1} A^{-1} $, microlocally near $ ( 0 , 0)$. 
That is the same as solving
\[  [  h D_{x_1}, A  ] + E A = 0 \,, \]
Since the principal symbol of $ P $ is independent of $ h$, same 
is true for the principal symbol of $ E$, $e$. Hence we can find 
$ a \in S^{0,0} ( T^* \RR^n ) $, independent of $ h$, 
$ a ( 0 , 0 ) \neq  0 $, and such that 
\[ \frac{1}{i} \{ \xi_1 , a \} + e a = 0 \, \]
near $ ( 0 , 0 )$. Choosing $ A_0 $ with the principal symbol $ a$
we can now find $ A_j \in \Psi^{ -j , 0 }_h ( T^* \RR^n ) $ so that
\[  [  h D_{x_1}, A_0 + A_1 + \cdots A_N  ] + E (A_0 + A_1 + 
\cdots A_N ) \in \Psi^{-N, 0 }_h ( T^* \RR^n ) \, .\]
We then put $ A \sim A_1 + A_2 + \cdots + A_N + \cdots $ which is
elliptic near $ ( 0 , 0 )$, and finally $ T = A^{-1} T_0 $. 
\end{proof}
Using the proposition we can transplant objects related to  $ P$
to the much easier to study objects related to $ hD_{x_1} $. 
In particular, we can microlocally define 
\begin{equation}
\label{eq:ker}
 {\rm ker}_{m_0 } ( P ) \stackrel{\rm{def}}{=} T^{-1} (
{\rm ker} \; (hD_{x_1}) )  \,,  \ \
{\rm ker} \; (hD_{x_1}) ) = \{ u \in {\mathcal D}' ( \RR^n ) \; 
: \; h D_{ x_1 } u = 0 \} \,. 
\end{equation}
Since $ {\rm ker} \; (hD_{x_1}) ) $ can be identified with 
$ {\mathcal D}' ( \RR^{n-1} ) $ we can also identify $  {\rm ker}_{m_0 } ( P )$
with $ {\mathcal D}' ( \RR^{n-1} ) $, 
microlocally near $ ( m_0 , ( 0 , 0 )) $:
\begin{equation}
\label{eq:K}
 K \; : \; {\mathcal D}' ( \RR^{n-1} ) \; \longrightarrow 
 {\rm ker}_{m_0 } ( P )  \,, \ \ K = T^{-1} \pi^* \,, \ \ 
\pi \; : \; x \mapsto ( x_2, \cdots x_n ) \,. \end{equation}

\section{Quantum time and quantum monodromy}
\label{qt}

Let $ P ( z ) \in \CI ( I_z ; \Psi^{0 , m } _h ( X ) ) $,
$ I = ( -a , a ) \subset \RR $,  be a smooth family of real 
principal type operators, with principal symbols, $ p ( z ) $, 
independent of $ h$.  We will assume that
\begin{gather*}
  \Sigma_{p(z)} \stackrel{\rm{def}}{=}
\{ m \in T^* X \; : \; p (z)( m ) = 0 \} \Subset T^* X \,, \ 
\text{ for $ z \in I $} \\
\text{ $ P ( z) $ is formally self-adjoint for $ z \in I $.} 
\end{gather*}
We assume that $ m_0 ( z) $ is a smooth 
family of periodic points of $ H_{p(z)} $,
with the minimal periods $ T( z)$ also smooth in $ z$, 
and the orbits $ \gamma (z )$:
\[  \exp ( T (z) H_{p(z)} ) ( m_0 (z)  ) = m_0 (z) \,, \ \ 
 \gamma \stackrel{\rm{def}}{=} \{ \exp ( t H_{p(z)} ) ( m_0 (z)  ) \; 
: \;  0 \leq t \leq T (z) \} \,. \]
When no confusion is likely to arise we may drop the dependence on $z$
in the notation.

 Let $ \Omega $ be a neighbourhood of $ \gamma (0)$
in $ T^* X $, 
\[ \Omega \simeq \gamma (0) \times {\mathbb B}_{ \RR^{2n-1} } ( 0 , 
\epsilon ) \,,\]
and we assume that for $ z \in I $, the orbits $ \gamma ( z) $ are
also contained in $ \Omega$. 
We now introduce a covering space of this tubular neighbourhood
\[ \widetilde \Omega \simeq \RR \times  {\mathbb B}_{ \RR^{2n-1} } ( 0 , 
\epsilon ) \,, \ \ \pi : \widetilde \Omega \longrightarrow \Omega \,,\]
with the lift of $ p ( z ) $ denoted by $ \widetilde p ( z )$, and we
will use the same notation for other objects. 

We start with the following
\begin{lem}
\label{l:3.1}
The tubular neighbourhood, $ \Omega $, of $ \gamma(0)$,  can be chosen 
small enough, so that the cover $ \widetilde \Omega $ contains no closed
orbits of $H_{\widetilde p ( z ) } $, $ z \in [ - \delta, \delta ] \subset
I $, for some small $ \delta > 0 $. 
\end{lem}
\begin{proof}
Let $ m \mapsto \tilde t ( m ) $ be a smooth function on $ \widetilde 
\Omega $ with the property that $ \tilde t ( \exp( t H_{ \widetilde p ( 0)} )
) = t $, and that $ d \tilde = \pi^* d \hat t  $, where $ d \hat t $ is a 
well defined one form in $ \Omega $. Then $ H_{ \widetilde p ( 0 ) } \tilde
t > 0 $ on the lift of $ \gamma $, and by shrinking $ \widetilde \Omega $
if necessary we conclude that $ H_{\widetilde p ( 0 ) } \tilde t > 0 $ 
on $ \widetilde \Omega $. By the periodicity and and a compactness 
argument we conclude that this holds for $ 0 $ replaced by $ z \in 
[ - \delta , \delta]$. Hence there are no closed orbits of 
$ H_{\widetilde p ( z ) } $ in $ \widetilde \Omega $. 
\end{proof}

We will now replace  $ \widetilde \Omega $ by 
a finite part:  $ \widetilde \Omega \simeq [ - L , L ] 
\times  {\mathbb B}_{ \RR^{2n-1} } ( 0 , \epsilon ) $, $ L \gg T $.

A {\em classical time function}, 
$ \tilde q ( z ) \in \CI ( \widetilde \Omega ; \RR)$, 
on $ \widetilde \Omega $ is defined as a solution of 
\begin{equation}
\label{eq:3.1}
\partial_z \tilde p ( z ) = - \{ \tilde p (  z ) , \tilde q ( z ) \} \,.
\end{equation}
In view of Lemma \ref{l:3.1} this equation can be solved (strictly 
speaking that may involve shrinking $ \widetilde \Omega $ further depending 
on the initial data,  but 
for simplicity of exposition we will ignore this point), and we can 
in particular consider solutions satisfying $ \widetilde q ( \tilde 
m_0 ( z ) ) 
= 0 $. In a neighbourhood of $ m_0 = m_0 ( 0 ) \in \Omega $ we can 
define $ q ( z) \in \CI $ such that
\[ \tilde q ( z) = \pi^* q ( z) \,, \
\text{ near $ \tilde m_0 $}, \ \   q (z) ( m_0(z) ) = 0  
 \ \ \pi : \widetilde 
\Omega \rightarrow \Omega \,. \]
We clearly have $ \partial_z p = \{ p , q \} $ near $ m_0 $. 
This defines the {\em local classical time} near $ m_0 $.
We also define the {\em first return classical time} near $ m_0$ by
demanding that 
\[ \tilde q ( z) = \pi^* q_{\circlearrowright}
 ( z) \,, \  
\text{ near $ \exp T ( 0 ) H_{ \tilde p ( 0 ) } 
\tilde m_0 $,} \ \
q_{\circlearrowright} (z) ( m_0(z) ) = \tilde q 
( \exp T( z ) H_{ \tilde p ( z ) } ( \tilde m_0 ( z )  ) \,.\]
An iteration procedure similar to the one recalled in the proof
of Proposition \ref{p:2.1} gives the quantum analogues microlocally 
defined near $ m_0$:
\begin{eqnarray}
\partial _z P ( z ) & = & - \frac{i}{h} [ P ( z ) , Q ( z) ] \,, 
\ \ \sigma ( Q ( z ) ) = q ( z ) \,, 
\label{eq:3.3}
\\
\partial _z P ( z ) & = & - \frac{i}{h} [ P ( z ) , Q_{\circlearrowright}
 ( z) ] \,, 
\ \ \sigma ( Q_{\circlearrowright} ( z ) ) = q_{\circlearrowright} ( z ) \,, 
\label{eq:3.3'}
\end{eqnarray}
Replacing $ Q ( z) $ by $ ( Q ( z) + Q ( z) ^*)/2 $, we can assume that
$ Q( z ) $ is formally self-adjoint. We clearly have 
\[   Q_{\circlearrowright} ( z) - Q ( z)  \; : \; {\rm{ker}}_{m_0} P ( z) 
\ \longrightarrow \  {\rm{ker}} 
_{m_0} P ( z)  \,. \]
Then $ Q ( z) $ is the {\em quantum time } near $ m_0 $, and 
$  Q_{\circlearrowright} $ is the {\em first return quantum time} near
$ m_0 $. See the proof of Lemma \ref{l:6.2} for further discussion of these
objects in the classical context.

For $ ( z, w ) $ near $ ( 0 , 0 ) $, and microlocally near $ m_0 $,
we can solve the following system of equations
\begin{equation}
\label{eq:3.4}
\begin{split}
& ( h D_z - Q ( z ) _L ) U ( z, w ) \stackrel{\rm{def}}{=} 
h D_z U ( z, w ) - Q ( z ) U ( z , w ) = 0  \\
& ( h D_w + Q ( w ) _R ) U ( z, w ) \stackrel{\rm{def}}{=} 
h D_z U ( z, w ) +  U ( z , w )  Q ( w )= 0  
\end{split}
\end{equation}
with the initial condition $ U ( 0 , 0 ) = Id $, and with $ 
U ( z , w ) $ bounded on $ L^2 $ (microlocally near $ ( m_0 , m_0 )$):
the solvability of the system follows from the fact that
\[ [ h D_z - Q( z)_L , h D_w + Q( w ) _R ] = 0 \,. \]
We easily check that (as always, microlocally)
\begin{equation}
\label{eq:3.5}
 U ( z , z ) = Id \,, \ \ U ( z, w ) U ( w , v ) = U ( z, v) \,, 
\end{equation}
and that $ U ( z, w ) $ is unitary. 
In fact, 
\[ h D_z ( U ( z , z )) =  Q ( z ) U ( z , z ) - U ( z , z ) Q ( z ) 
= - [ U ( z , z ) , Q ( z ) ] \,, \ \ U ( 0 , 0 ) = Id \,, \]
and $ U ( z , z ) = Id $ is the unique solution.  The other property 
is derived similarly:
\begin{gather*}
 h D_w ( U ( z , w ) U ( w, v ) ) = - U ( z , w ) Q  (w ) U ( w, v ) 
+ U ( z , w ) Q ( w ) U ( w, v ) = 0 \,, \\ 
U ( z , w ) U ( w, v  ) 
\rest_{ w = z } = U ( z , v ) \,. \end{gather*}

By varying $ m_0 $ along the orbit of $ H_{p(0)}$, and 
by extending $ Q ( z ) $ maximally forward ($+$) and backward ($-$),
we can define semi-global versions of $ U ( z , w ) $:
\[ \begin{split}
& U_+ ( z, w )  \ \text{ microlocally on a neighbourhood of the diagonal over} 
\\
& \ \ \ \ \ \ \ \ \ 
\{ \exp t H_p ( m_0 ) \; : \; 
 - \epsilon <  t < T(0)  - 2 \epsilon  \}
  \,, \\
& U_- ( z, w )  \ \text{ microlocally on
a neighbourhood of the diagonal over
} 
\\
& \ \ \ \ \ \ \ \ \ 
\{ \exp t H_p ( m_0 ) \; : \; 
 - \epsilon < - t < T(0)  - 2 \epsilon  \}
 \,.
\end{split} 
\]

The operators have the following intertwining property:
\begin{prop}
\label{p:3.1}
Microlocally near the diagonal over 
\[ \{ \exp t H_p ( m_0 ) \; : \; 
 - \epsilon < \pm  t < T(0)  - 2 \epsilon  \} \,, \]
 and for $z,w $ close to $ 0$, we have 
\[ P ( z) U_\pm  ( z , w ) = U_\pm  ( z, w ) P ( w ) \,.\]
\end{prop}
\begin{proof}
We define $ P^\sharp ( w ) = U ( w, z ) P ( z) U ( z , w ) $ and 
differentiate with respect to $ w $:
\begin{gather*}
h D_w P^ \sharp ( w) =  Q ( w ) U ( w, z) P ( z ) U ( z ,w ) - 
U ( w, z ) P ( z ) U ( z, w ) Q ( w ) = 
- [ P^\sharp ( w ) , Q ( w ) ] \,, \\ P^\sharp ( w ) \rest_{ w = z} = 
P ( z ) \,, \end{gather*}
that is, $ P^\sharp ( w ) $ satisfies \eqref{eq:3.3} and consequently 
$ P^\sharp ( w ) = P ( w ) $. 
\end{proof}

By replacing the local quantum time, $ Q ( z ) $, 
 by the first return quantum time, $ \QC ( z ) $ (see 
\eqref{eq:3.3}, \eqref{eq:3.3'}), we also define $ \UC ( z, w) $,
\[ \UC ( z ,z ) = Id\,, \ \ \UC ( z , w ) \UC ( w, v ) = 
\UC ( z , v ) \ \text{ microlocally near $ m_0 $.}\]
This definition will be useful when we study the {\em quantum monodromy 
operator}. To introduce it, we first define the forward and backward
propagators:
\begin{gather}
\label{eq:epm}
\begin{gathered}
 I_\pm ( z ) \; : \;  {\rm ker}_{m_0( z )  } ( P( z )  ) 
\ \longrightarrow \ {\mathcal D}' ( X ) \\
I_\pm ( z ) = Id
\rest_{{\rm ker}_{m_0( z )  } ( P( z )  )}
\,,  \ \text{ microlocally near $ m_0 ( z )$,} \\
P ( z ) I_\pm ( z ) = 0 \,,  \ \text{ microlocally near $
\{ \exp ( tH_{ p ( z ) } ) m_0 ( z ) \; : \;  -\epsilon < \pm t < T ( z ) - 
2 \epsilon  \} $.} 
\end{gathered}
\end{gather}
That the operators $ I_\pm ( z ) $ are microlocally well defined follows from 
Proposition \ref{p:2.1}, and ``microlocally'' is meant via the 
identification of 
$ {{\rm ker}_{m_0( z )  } ( P( z )  )} $ with $ {\mathcal D' } ( 
\RR^{n-1}) $ as in \eqref{eq:K}. 
We fix $ m_0 = m_0 ( 0 )$ as in the definition of $ \widetilde 
\Omega $ above, and define
\begin{equation}
\label{eq:wpm}
\begin{split}
& W_+ = \ \text{ a neighbourhood of $ m_0 $ in $ T^* X $,}\\
& W_- = \ \text{ a neighbourhood of $ \exp ((T(0)/2) H_{ p (0) } ) m_0 $ 
in $ T^* X $,}\\
& W_- \subset \bigcup_{ | t + T(0)/2 | < \epsilon } \exp t H_{ p ( 0) } 
( W_+) \,,
\end{split}
\end{equation}
noting that for $ z $ small enough, we can replace $ m_0$, $ T (0) $,
$ p ( 0  ) $, by  $ m_0( z ) $, $ T (z ) $, $ p ( z  ) $ in this definition.
This shows that $ I_- ( z ) $ maps $ {\rm{ker}}_{m_0 ( z ) } ( P ( z ) )$
onto $ {\rm{ker}}_{ \exp ( (T(z ) /2 ) H_{ p ( z)} ) ( m_0 ( z ) ) } 
( P ( z ) )$, microlocally near $ W_- \times W_+ $. This means that the 
left microlocal inverse exists 
and we can give the following

\medskip

\noindent
{\bf Definition.} {\em
The (absolute) quantum monodromy operator
\[   {\mathcal M} ( z)  \; : \;   {\rm{ker}}_{m_0 ( z ) } ( P ( z ) )
\longrightarrow   {\rm{ker}}_{m_0 ( z ) } ( P ( z ) ) \,, \]
is microlocally defined near $ W_+$ by
\begin{equation}
\label{eq:mon}
 I_+ ( z ) f = I_- ( z) {\mathcal M} ( z ) f \,, \ \ 
f \in   {\rm{ker}}_{m_0 ( z ) } ( P ( z ) )\,,
 \ \text{ microlocally near $ W_-$.}
\end{equation}
The quantum monodromy operator,
\[ M ( z ) \; : \; {\mathcal D}' ( \RR^{n-1} ) \rightarrow 
{\mathcal D}' ( \RR^{n-1} ) \,, 
\] 
is microlocally defined near $ ( 0 , 0 ) \in T^* \RR^{n-1} $, by 
\begin{equation}
\label{eq:3.mon}
M ( z ) =   K( z ) ^{-1} {\mathcal M } (z ) K ( z) \,,\end{equation}
where $ K ( z) $ is as in \eqref{eq:K}.
}

\medskip

The basic properties are given in 
\begin{prop}
\label{p:3.m}
Let $ U ( z , w ) $ and $ \UC ( z, w ) $ be given by \eqref{eq:3.3},
\eqref{eq:3.3'}, and \eqref{eq:3.4}. Then the following diagram commutes
(microlocally near $ m_0$):
\begin{equation}
\label{eq:3.6}
\begin{array}{ccc}
{\rm{ker}}_{m_0 ( w ) } ( P ( w ) ) & 
\stackrel{{\mathcal M}(w)}{\longrightarrow} & 
{\rm{ker}}_{m_0 ( w ) } ( P ( w ) )\\
\Big\downarrow\vcenter{%
\rlap{$U ( z, w)$}}& & 
\Big\downarrow\vcenter{%
\rlap{$\UC ( z, w)$}} \\
{\rm{ker}}_{m_0 ( z ) } ( P ( z) ) 
 & 
\stackrel{{\mathcal M}(z)}{\longrightarrow} & 
{\rm{ker}}_{m_0 ( z ) } ( P ( z ) )
\end{array}
\end{equation}
Choosing $ K ( z ) $, so that $ K ( z ) = U ( z, w ) K( w) $, we also have 
\begin{equation}
\label{eq:3.m}
 h D_z M ( z ) = [ K ( z )^{-1} ( \QC ( z ) - Q ( z ) ) K ( z ) ]
M ( z) \, ,\end{equation}
where we recall that by \eqref{eq:3.3} and \eqref{eq:3.3'}, 
$ \QC ( z ) - Q ( z ) \; : \;  {\rm{ker}}_{\tilde m_0 (z ) } ( 
\widetilde P ( z ) ) \rightarrow  {\rm{ker}}_{\tilde m_0 ( z ) } ( 
\widetilde P ( z ) ) \,, $ and hence $  K ( z) ^{-1} $ is well 
defined.

\end{prop}
\begin{proof}
We need to show that $ \UC ( z , w )  \MO ( w ) = \MO ( z ) 
U ( z , w ) $, and since $ \UC $ is naturally defined using the 
covering space, we will translate this into a statement there. 
We can microlocally define $ \widetilde P ( w ) $ on $ \widetilde
\Omega $ and then,
\[ \widetilde I_+ ( w)
 \; : \; {\rm{ker}}_{\tilde m_0 ( w ) } ( 
\widetilde P ( w ) ) \longrightarrow {\rm{ker}} ( \widetilde P ( w ) )\,,\]
and we define
\[ \widetilde \MO ( w )  \; : \; 
 {\rm{ker}}_{\tilde m_0 ( w ) } ( 
\widetilde P ( w ) ) \longrightarrow
 {\rm{ker}}_{\exp ( T ( w ) H_{ \tilde p ( w )} ) (\tilde m_0 ( w )) } ( 
\widetilde P ( w ) ) \,, \]
by restricting $ \widetilde I_+ ( w ) $ to a neighbourhood of $ 
\exp ( T ( w ) H_{ \tilde p ( w )} ) (\tilde m_0 ( w ))$. Since 
for $ \pi : \widetilde \Omega \rightarrow \Omega $, we microlocally have
\begin{gather*}
 \pi_* \; : \;  {\rm{ker}}_{\exp ( T ( w ) H_{ \tilde p ( w )} ) (\tilde m_0 ( w )) } ( 
\widetilde P ( w ) ) \longrightarrow {\rm{ker}}_{m_0 ( w ) } P ( w ) \,, \\
 \pi_* \widetilde \MO (w) \pi^* = \MO ( w) \,.
\end{gather*}
Using the quantized version of $ \tilde q $ in \eqref{eq:3.1}, we 
also define $ \widetilde U ( z, w)$, so that 
$ \widetilde U ( z, w ) \widetilde P ( w ) = \widetilde P ( z ) 
\widetilde U ( z , w) $. In particular we have 
\[ \widetilde U ( z, w ) \widetilde I_+ ( w ) = \widetilde I_+ ( z ) 
\widetilde U ( z , w) \,. \]
Restricting (microlocally) to a neighbourhood of 
\[ ( \exp ( T ( z ) H_{ \tilde p ( z )} ) (\tilde m_0 ( z )) , 
\tilde m_0 ( z ) ) \in \widetilde \Omega \times \widetilde \Omega \,, \]
and projecting to $ \Omega \times \Omega $, we obtain
\[ \UC ( z, w ) \MO ( w ) = \MO ( z ) U ( z, w ) \,.\]

To see \eqref{eq:3.m} we first note that differentiation of $ K ( z ) 
= U ( z, w ) K ( w) $ and the definition of $ U ( z , w ) $ gives
\[ h D_z K ( z ) =  Q ( z) K ( z ) \,. \]
We then use the commutative diagram to see that
\[ K ( z ) M ( z) = \UC ( z , w ) {\mathcal M} ( w ) U ( w, z ) K ( z ) \,.\]
Differentiating this with respect to $ z $ and using the previous equation
gives
\[  K ( z) h D_z M ( z) = 
 ( \QC ( z ) - Q ( z ) ) K ( z ) M ( z) \, . \]
We then recall that by \eqref{eq:3.3} and \eqref{eq:3.3'}, 
$ Q ( z ) - \QC ( z ) \; : \;  {\rm{ker}}_{\tilde m_0 (z ) } ( 
\widetilde P ( z ) ) \rightarrow  {\rm{ker}}_{\tilde m_0 ( z ) } ( 
\widetilde P ( z ) ) \,, $ and hence $  K ( z) ^{-1} $ can be 
applied to both sides.
\end{proof}

We can define the {\em Poincar{\'e} map} for $ \gamma $ with primitive 
period $ T$:
\[ C : T^* \RR^{n-1} \longrightarrow T^* \RR^{n-1} \,, \ \
\text{defined near $ ( 0 , 0 )$}\,, \ \ C ( 0 , 0 ) = ( 0, 0 ) \,,
\]
as follows: for a neighborhood of  $ m_0 \in \gamma  $, $ U_0$, 
$ U_0 / \exp({t H_{ p} }) $ can be identified with a neighbourhood
of $ ( 0 , 0 ) \in T^* \RR^{n-1} $ (using the local identification 
of $ p $ with $ \xi_1 $, as in the proof of Proposition \ref{p:2.1}),
with $ [m_0] $ corresponding to $ ( 0 , 0 )$. The Poincar{\'e} map
is then given by 
\begin{gather}
\label{eq:4.poin}
\begin{gathered}
C : \kappa^{-1} ([m])  \longmapsto  \kappa^{-1} ( [ \exp ( T H_p ) m ] ) \,, 
 \\ [m] \in U_0 / \exp ( t H_p ) \,, \kappa : T^* \RR^{n-1} \rightarrow
U_0 / \exp ( t H_p ) \,. \end{gathered} \end{gather}
It will always be undestood that $ \kappa $ chosen here is the 
symplectic transfomation corresponding to $ K = K( z) $ in \eqref{eq:K}
and \eqref{eq:3.mon}.

To study quantum properties of the monodromy operator
it is convenient to  introduce $ \chi \in \CIc ( T^* X ) $ satisfying
\begin{gather}
\label{eq:4.chi}
\begin{gathered}
\chi \equiv \left\{ \begin{array}{ll}
1 &  \text{near $\{ \exp( t H_{p(0)} (m_0) ) : \epsilon < t < T(0)/2 - 
\epsilon \}$} \\
0 &  \text{near $\{ \exp( t H_{p(0)} (m_0) ) : \epsilon < -t < T(0)/2 - 
\epsilon \}$} \end{array} \right. \\
\Omega \cap \{ m \; : \; \chi ( m ) \neq 1 \} 
\cap \{ m \; : \; \chi ( m ) \neq 0 \} \subset  W_+ \cup W_- \,, 
\end{gathered}
\end{gather}
where $ W_\pm $ are as in \eqref{eq:wpm}, and $ \Omega $ is a
small neighbourhood of $ \gamma $.  If $ \rho_\pm \equiv 1 $
microlocally near $ W_\pm $, and $ \rho_\pm \equiv 0 $ near 
$ W_\mp $, we define
\[ [ P, \chi]_{W_\pm} = \rho_\pm  [ P, \chi ] \,,\]
where we use the same notation for $ \chi $ and $ {\rm{Op}}_h( \chi )$. 
We then have the basic property of the {\em quantum flux} (see 
\cite{Gr}):
\begin{lem}
\label{l:4.1}
Let $ K ( z) $ be in \eqref{eq:K}. Then \[ U ( z ) \stackrel{\rm{def}}{=} 
K( z) ^* [ ( i/h) P ( z) ,\chi]_{W_+} K( z) \; : \; {\mathcal D}'( \RR^{n-1})
\; \longrightarrow \;  {\mathcal D}'( \RR^{n-1})   \]
is microlocally positive near $ ( 0 , 0 ) \in T^* \RR^{n-1} $ and independent
of $ \chi $ with the properties \eqref{eq:4.chi}.

If we replace $ K ( z) $ by $ K( z) U(z)^{-\frac12} $ then 
\begin{equation}
\label{eq:4.2}
K( z) ^* [ ( i /h ) P (z) ,\chi ]_{W_+}  K ( z) = Id 
\ \text{ microlocally near $ ( 0 , 0 ) \in T^* \RR^{n-1} $.}
\end{equation}
\end{lem}
\begin{proof}
We note that if $ P ( z) u = 0 $ near $ W_+ $, and $ \tilde \chi$ is 
another function satisfying \eqref{eq:4.chi}, then 
\[ K( z)^* [  ( i /h ) P ( z) , \chi - \tilde \chi  ] u = 
K ( z )^* P ( z) ( \chi - \tilde \chi ) u - K(z)^* ( \chi - \tilde \chi) 
P ( z ) u = 0 \,,\]
since $ P ( z) u = 0 $, and $ K ( z ) ^* P ( z ) 
=  ( P ( z) K( z))^* = 0 $. The positivity also comes from expanding
the commutator and using Proposition \ref{p:2.1}:
\[ \langle \tilde K ( z)^* [ (i/h) hD_{x_1} , \chi] _{W_+} \tilde K ( z)  u , 
 u \rangle = \langle \partial_{x_1} \chi \rho_+  \tilde K ( z ) u , 
\tilde K ( z ) u \rangle \geq \langle \tilde \rho u , \tilde \rho u \rangle
\,, \]
where $ \tilde K ( z) $ is the composition of $ K ( z) $ and  $ T$ of 
Proposition \ref{p:2.1}, and $ \tilde \rho \equiv 1 $ in a neighbourhood
of $ ( 0, 0) \in T^* \RR^{n-1} $ (we again use the same notation for the
function and its quantization). 
\end{proof}

From now on, our choice of $ K ( z) $ in \eqref{eq:K} is made so that
\eqref{eq:4.2} holds.  We only need to check that we still have 
\[ K ( z ) =  U ( z, w ) K  ( w ) \,. \]
In fact, we have in general, in the microlocal sense,
\[ \begin{split}
& K ( z )^* [ ( i / h ) P ( z) , \chi ]_{ W_+} K ( z ) = \\
& K ( w )^* U ( w , z ) [ ( i / h ) P ( z) , \chi ]_{W_+} 
U ( z, w ) K ( w) = \\
& K( w)^* [ ( i /h ) P ( w ) , \tilde \chi ]_{ W_+} K ( w )\,, 
\end{split} \]
and the last expression is unchanged if we replace $ \tilde \chi $ 
by $ \chi$ (the quantum flux property used before). We also used the
unitarity of $ U ( z, w ) $. 

With this choice of $ K ( z) $ we have the following important
and well known
\begin{prop}
\label{p:3.2}
The monodromy operator, $ M ( z )$, defined by \eqref{eq:3.mon} with 
$ K ( z) $ satisfying \eqref{eq:4.2} is microlocally unitary:
\[ M ( z) ^* = M( z )^{-1} \ \ \text{ microlocally near $ ( 0, 0 ) \in 
T^* \RR^{n-1} $} \,, \]
and it is an $h$-Fourier Integral Operator:
\[ M ( z ) \in I^0 ( \RR^{n-1} \times \RR^{n-1} ; {C}( z) ')\,,\]
where $ { C} ( z) $ is the {\em Poincar{\'e} map} \eqref{eq:4.poin}.
\end{prop}
\begin{proof}
We need to show that for $ v \in {\mathcal D}' ( \RR^{n-1} ) $ with 
$ WF_h ( v ) $ in a neighbourhood of $ ( 0, 0 )$, we have 
\begin{equation}
\label{eq:3.7} \langle M ( z) v , M ( z ) v \rangle = \langle v , v \rangle + 
{\mathcal O} ( h^\infty ) \| v \|^2 \,. \end{equation}
If we put $ u = K ( z ) v $, use \eqref{eq:4.2}, and the defintion of 
$ M ( z ) $, \eqref{eq:3.mon}, then the
left hand side of \eqref{eq:3.7} becomes:
\[ \langle K( z )^* [ P ( z)  , \chi]_{W_+} {\mathcal M } ( z) u , 
K ( z ) ^{-1} {\mathcal M } ( z ) u \rangle  = 
\langle [ P ( z)  , \chi]_{W_+} {\mathcal M } ( z) u , 
 {\mathcal M } ( z ) u \rangle
,.\]
As in the proof of Lemma \ref{l:4.1}, we see that 
for 
$ 0 < t < T ( 0 )/2  + \epsilon $, the right hand side of the previous 
expression is equal to (modulo $ {\mathcal O} ( h^\infty ) $)
\[ \langle  [ ( i / h ) P ( z) , (\exp ( tH_{p(0)} )^* \chi]_{
\exp (- t H_{ p( 0 )} ) W_+} I_ - ( z ) {\mathcal M } ( z) u 
, I_- ( z ) {\mathcal M } ( z) u  \rangle 
\,, \]
which corresponds to moving the support of $ \chi$ in the direction opposite
to the flow of $ H_{p ( 0 )} $, and simultaneously moving $ W_+ $ so that
\eqref{eq:4.chi} holds.

Similarly, for $ - T ( 0 )/2  - \epsilon < t < 0 $, the right hand side of 
\eqref{eq:3.7} is equal to 
\[ \langle [ ( i / h ) P ( z) , (\exp (- t H_{p(0)} )^* \chi]_{
\exp (   t   H_{ p( 0 )} ) W_+} I_+ ( z) u ,  I_+ ( z ) \rangle \,. \] 
For $ t \sim T ( 0)/ 2 $, $ \exp ( \pm t H_{ p ( 0 ) } ( W_+ ) \subset
W_- $,
and hence 
\begin{gather*}
 \langle [ ( i / h ) P ( z) , (\exp (- t H_{p(0)} )^* \chi]_{
\exp (   t   H_{ p( 0 )} ) W_+} I_+ ( z) u ,  I_+ ( z ) u \rangle =\\
 \langle  [ ( i / h ) P ( z) , (\exp ( tH_{p(0)} )^* \chi]_{
\exp (- t H_{ p( 0 )} ) W_+} U_ - ( z ) {\mathcal M } ( z) u 
, I_- ( z ) {\mathcal M } ( z) u  \rangle \,, \ \ t \sim T ( 0)/ 2\,,
\end{gather*} 
from the definition of $ {\mathcal M}( z ) $, \eqref{eq:mon}.
But this shows \eqref{eq:3.7} proving the first part of the proposition.

To see the second part we use use Lemma \ref{l:egor}, and the obvious
conjugation properties of the solution in the model case discussed in 
Proposition \ref{p:2.1}: going around the closed orbit we obtain that 
the underlying symplectomorphism is given by the Poincar{\'e} map. 
\end{proof}

So far we have discussed only the case of $ z \in \RR $. We can 
now consider almost analytic extensions of the operators 
$ Q ( z ) $, $ \QC ( z ) $, $ U_\pm (z, w )$, $ I_\pm ( z ) $,
and $ M ( z )$. For that we consider a complex neighbourhod of
$ I \subset \RR $:
\[ I_{ h , L} = \{ z \; : \; \Re z \in I \,, \ |\Im z | \leq 
L h \log ( 1/h ) \} \,.\]
The families of pseudo-differential operators $ P ( z) $, $ Q ( z )$, 
and $ \QC ( z ) $ have almost analytic extensions given by Lemma \ref{l:2.1},
and we use the same notation for them.  We then use Lemma \ref{l:2.2} 
and \eqref{eq:3.4} to 
extend $ U ( z,w ) $, $ \UC ( z , w )$, and $ U_\pm ( z , w ) $ to 
$ ( z, w ) \in I_{ h , L } \times I_{ h , L } $. We then have
\[ P ( z ) U _\bullet( z ,w ) = U_\bullet ( z, w ) P ( w) \,, \ \
 ( z, w ) \in I_{ h , L } \times I_{ h , L } \,, \]
microlocally (that is, in particular modulo $ {\mathcal O} ( h^\infty ) $). 
Indeed, for $ x, w \in I\times I $, and $ |y | \leq L h \log ( 1 /h ) $, 
we have, as in the proof of Lemma \ref{l:2.2},
\begin{gather*}
\partial_y [ P ( x + i y ) U ( x + i y, w ) - U ( x + i y , w) P ( w ) ]
= {\mathcal O} ( y^\infty ) \,, 
\\
 [ P ( x + i y ) U ( x + i y, w ) - U ( x + i y ) P ( w ) ]\rest_{ y = 0 }
= 0 \,. \end{gather*}
Hence we can define
\[ I_\pm ( z ) = U_\pm ( z , w) I _\pm ( w ) \,, \ 
( z , w ) \in I_{h, L } \times I \,, \]
so that $ P ( z ) I  _\pm ( z ) = 0 $. 

To define an almost analytic extension of $ M ( z )$ we first almost 
analytically extend the pseudo-differential operator 
$ K ( z) ^{-1} ( Q  ( z ) - \QC ( z) ) K ( z ) $, and then use 
\eqref{eq:3.m} and Lemma \ref{eq:2.2}. In particular, Proposition 
\ref{p:3.2} gives,
\[ M ( z) ^{-1} = M ( \bar z )^* \,, \ \  |\Im z | \leq L h \log ( 1/h) \,.
\]

\section{Grushin problem near a closed trajectory}
\label{gr}

As in the previous section 
we assume that $ P ( z )$ is {\em self-adjoint} for $ z \in \RR $, and
denote by the same symbol the almost analytic continuation of $ P ( z )$.
Although the inverse of $ P ( z) $ does not 
normally exist near $ \gamma = \gamma ( 0 ) $ for all $ z \in I $ we will
describe $ P ( z )^{-1} $ in terms of the inverse of a microlocal 
effective Hamiltonian $ E_{- +} ( z ) = I - M ( z )$. We will do it
first for $ z $ real and then use the extensions of operators 
$ U _\pm ( z , w )$ described at the end of the last section to 
transplant the results to complex values of $ z$.

To do that we follow the now standard Grushin reduction \cite{Gr},
and consider the system
\begin{equation}
\label{eq:4.p} {\mathcal P} ( z ) = \left( \begin{array}{ll} 
 ( i / h ) P ( z) & R_- ( z ) \\
\; \ R_+ ( z ) & \ \ 0 \end{array} \right) \; : \; 
{\mathcal D}' (  X ) \times {\mathcal D}' ( \RR^{n-1}) 
\; \longrightarrow \;
{\mathcal D}' (  X ) \times {\mathcal D}' ( \RR^{n-1}) \,,\end{equation}
defined microlocally near $ \gamma \times {(0, 0)} $, and where the
operators $ R_\pm $ need to be suitably chosen. 

We will successively build the operator $ {\mathcal P}(z)$ and its
inverse. We start by putting
\begin{equation}
\label{eq:4.3}
 R_+ ( z) = K( z )^* [ ( i/h ) P  ( z) , \chi ]_{W_+} \,,\end{equation}
and $ u $ with $ P u = 0 $ near $ W_+$, $ R_+ (z)  u $, is its 
Cauchy data. Hence $  u = K ( z) v $ provides a {\em local}
solution to the microlocal Cauchy problem:
\begin{equation}
\label{eq:4.4}
\left\{ \begin{array}{l} P ( z ) u = 0 \\
R_+ ( z ) u = v \end{array} \right. 
\end{equation}
To obtain a global Cauchy problem we need to introduce $ R_- ( z )$.
To do that we define
\[ K_f ( z ) = I_+ ( z ) K ( z )  \,, \ \ K_b ( z ) = I_- ( z ) K ( z ) \,,
\] 
where the operators $ I_\pm ( z) $ are defined in \eqref{eq:epm}. 
We recall the definition of the monodromy operator:
\begin{equation}
\label{eq:4.6}
  K_f  ( z) = K_b ( z ) M (z )   \text{ microlocally near $W_- \times {(0,0)}
 $.}
\end{equation}
We can build a solution of \eqref{eq:4.4}
 in $ \Omega \setminus W_-  $ by 
putting 
\begin{equation}
\label{eq:4.7}
E_+ (z ) v = \chi K_f (z )v  + ( 1 -\chi ) K_b (z) v \,,
\end{equation}
so that in particular, $ E_+ ( z ) v = K  ( z) v $, in $ W_+ $, and
consequently
\begin{equation}
\label{eq:4.9}
R_+ ( z ) E_+ ( z ) = Id  \text{ microlocally near $ ( 0, 0 ) \in 
T^* \RR^{n-1} $. } 
\end{equation}
Applying the operator, and using \eqref{eq:4.6} we obtain
\[ \begin{split} \frac{i}{h} P( z) E_+( z) v &  =  
[ ( i /h ) P ( z) , \chi]_{W_-} K_f ( z) v - 
[ ( i /h ) P ( z) , \chi]_{W_-} K_b ( z) v  \\
& = 
[ ( i /h ) P ( z) , \chi]_{W_-} K_b ( z) ( M (z) - I ) v \,. 
\end{split}
\]
Hence we obtain a globally (near $ \gamma$) solvable Cauchy problem 
by putting 
\begin{equation}
\label{eq:4.10}
\left\{ \begin{array}{l} \frac{i}h P ( z ) u + R_- ( z ) u_- = 0 \\
R_+ ( z ) u = v \end{array} \right. \,,
\end{equation}
with 
\begin{equation}
\label{eq:4.11}
R_- ( z) = [ ( i /  h ) P ( z) , \chi]_{W_- } K_b ( z ) \,. \end{equation}
The problem \eqref{eq:4.10} is solved by putting
\begin{equation}
\label{eq:4.12}
u = E_+ ( z )  v \,, \ \ u_- = E_{ -+ } (z) v \,, \ \ 
E_{ -+ } ( z) \stackrel{\rm{def}}{=} I - M ( z ) \,,
\end{equation}
where $ E_+ ( z ) $ was given by \eqref{eq:4.7}. 

The definitions
\eqref{eq:4.3} and \eqref{eq:4.11} give $ {\mathcal P} ( z ) $ in 
\eqref{eq:4.p}. If the microlocal inverse, $ {\mathcal E} ( z) $, exists,
it is necessarily given by
\begin{equation}
\label{eq:E}
 {\mathcal E} ( z) = \left( \begin{array}{ll}
E (z ) & E_+ ( z ) \\
E_- ( z) & E_{ - + } ( z) 
\end{array} \right)
\end{equation}
where $ E_+ ( z) $ and $ E_{-+} ( z) $ have already been constructed.

It remains to find $ E ( z)\,, E_- ( z)  $, 
and to show that the resulting operator
$ {\mathcal E} ( z ) $ is the right and left microlocal inverse of 
$ {\mathcal P } ( z ) $. For the right inverse, this means solving
\begin{equation}
\label{eq:4.18}
\left\{ \begin{array}{l} \frac{i}h P ( z ) u + R_- ( z ) u_- = v \\
R_+ ( z ) u = v_+ \end{array} \right. \,.
\end{equation}

We first introduce the forward and backward fundamental solutions of 
$ ( i / h ) P ( z) $:
\[ \begin{split}
& L_f ( z )  \ \text{ \ microlocally defined near  
$ ( \Omega \times_\epsilon \Omega )_{+}$ }\,, \\
& L_b ( z) \ \text{ \ microlocally defined near  
$ ( \Omega \times_\epsilon \Omega )_{-} $}  \,, \end{split} \]
where $ ( \Omega \times_\epsilon \Omega )_{\pm}$ is given by 
\begin{equation*}
  ( \Omega \times_\epsilon \Omega )_{\pm} \stackrel{\rm{def}}{=}
\left( \bigcup_{ m \in \Omega } \{ ( \exp ( t H_{p(0)}  )m , m )\right) 
\cap \Omega \times \Omega  \; : \; 
 - \epsilon < \pm t < T(0)  - 2 \epsilon  \}\,. \end{equation*}

 To do that we 
use Proposition \ref{p:2.1} and the corresponding local 
forward and backward fundamental solutions: 
\[ \begin{split}
& L_f^0 v ( x ) =  \int_{ -\infty }^{x_1} v ( t , x') dt \,,\\
& L_b^0 v ( x ) = -  \int_{x_1}^\infty  v ( t , x') dt \,,
\end{split}\]
$ v \in {\mathcal E}' ( \RR^{n} ) $.

We will now try to build an approximate  solution of $ ( i /h ) P ( z) u 
= v $ using $ L_\bullet ( z )$. For that let us put 
\[ \tilde u = L_f ( z ) ( 1  - \chi ) v \,.\]
Let us also define $ \chi_b $, $ \chi_f $ satisfying \eqref{eq:4.chi} and 
in addition,
\[ \text{ $ \chi_b = 1 $ on $ {\rm{supp}} \; \chi \cap W_+ $, 
$ \chi = 1 $ on $ {\rm{supp}} \; \chi_f \cap W_+ $. }\]
We now put 
\[\tilde u = L _f ( z ) ( 1 - \chi ) v \,, \]
where we can think of $ \tilde u $ as being microlocally defined on 
the covering space of $ \Omega $, $ \widetilde \Omega $ (see the proof of
Proposition \ref{p:3.m}). Hence, $  P( z ) \tilde u =  0 $ 
to the right of the support of $ 1 - \chi $ (in the direction of the 
flow), in particular on the support of $ \chi_f $. Hence, to the right of 
the support of $ 1 - \chi$, 
\[ \begin{split}
\tilde u & = K ( z) K( z) ^* [ ( i / h ) P ( z ) , \chi_f] _{W_+}
\tilde u \\
& = K ( z) K ( z)^* [( i /h ) P ( z) , \chi_f ] _{ W_+ } 
L_f ( z ) ( 1- \chi) v \,. \end{split}
\]
If we use the notation from the proof of Proposition \ref{p:3.m} and
put $ \widetilde K_f ( z ) = \widetilde I_+ ( z) K ( z ) $, then in the 
forward direction of propagation past the support of $ 1 - 
\chi $, we have in $\widetilde \Omega $, 
\begin{equation}
\label{eq:4.19}
\tilde u = \widetilde K_f ( z) K ( z) ^* [ ( i / h ) P ( z) , \chi_f ]_{
W_+} 
L_f  ( z) ( 1- \chi) v \,. \end{equation}
Similarly, if $ \hat u = L_b ( z) \chi v $ then, left to the support of 
$ \chi $, we have $ P ( z ) \hat u  = 0 $, and we can extend $ \hat u $ 
father left, microlocally in $ \widetilde \Omega $:
\begin{equation}
\label{eq:4.20}
\hat u = \widetilde K_b ( z)K ( z ) ^* [ ( i / h ) P ( z) , \chi_b ] 
_{W_+} L_b ( z) \chi v \,. \end{equation}
We can think of $ \tilde u $ and $ \hat u $ as multivalued in $\Omega $
and we will define, near $ W_- $,
\[ \begin{split}
& L_{ff} v = \text{ second branch of $ \tilde u $ near $ W_- $}\,,\\
& L_{bb} v = \text{ second branch of $ \hat u $ near $ W_- $}\,.
\end{split} \]
With this notation we put
\begin{equation}
\label{eq:4.21}
u_0 = E_0 ( z) v \stackrel{\rm{def}}{=} \left\{ 
\begin{array}{ll} L_b ( z ) \chi v + L_f ( z ) ( 1- \chi ) v & 
\text{outside $ W_-$  } \\
L_b ( z) \chi v + ( 1 - \chi ) L_{bb} ( z) \chi v + 
L_f ( z ) ( 1 - \chi ) v + \chi L_{ff} ( 1 - \chi ) v & 
\text{ in $ W_- $ } \end{array} \right. \,. 
\end{equation}
An application of $ ( i /h ) P ( z) $ gives
\[ ( i / h) P ( z) u_ 0 = v - [ ( i / h ) P ( z) , \chi ] _{W_-}
L_{bb} (z)  \chi v + [ ( i / h ) P ( z ) , \chi] _{W_-} L_{ff}( z)  
( 1 - \chi ) v \,,\]
and using \eqref{eq:4.19} and \eqref{eq:4.20} (where we now drop the
$ \tilde {}$ as we are taking the second branch of $ \tilde u $ and
$\hat u $,
 and the definition of $ M ( z) $,  we get
\[\begin{split}
\frac{i}{h} P ( z) u_0 & = v - [ ( i / h ) P(z ) , \chi ]_{W_-} 
K_b ( z) K( z )^* [ ( i / h ) P ( z) , \chi_b ] _{ W_+} L_b ( z) 
\chi v \\
& \ \ \ \ + [ ( i / h ) P(z ) , \chi ]_{W_-} 
K_f ( z) K( z )^* [ ( i / h ) P ( z) , \chi_f ] _{ W_+} L_f ( z) 
( 1 - \chi)  v  \\
& =  v - [ ( i / h ) P(z ) , \chi ]_{W_-} 
K_b ( z) \left( K( z )^* [ ( i / h ) P ( z) , \chi_b ] _{ W_+} L_b ( z) 
\chi v \right.\\ 
& \ \ \ \ - \left.
M( z)  K( z )^* [ ( i / h ) P ( z) , \chi_f ] _{ W_+} L_f ( z) 
( 1 - \chi)  v  \right)\,.
\end{split}
\]
In other terms,
\begin{equation}
\label{eq:4.23}
\frac{i}{h} P ( z ) E_0 ( z ) v + R_- ( z) E_{ 0, - } ( z) v = v \,,
\end{equation}
where we defined  $ E_0 ( z) $ by \eqref{eq:4.21} and 
\begin{equation}
\label{eq:4.24}
E_{0, -} ( z) = K( z) ^* [ ( i / h ) P ( z) , \chi_b ] L_b ( z) \chi 
- M ( z) K( z)^* [ ( i / h ) P ( z) , \chi_f ] _{ W_+} L_f ( z) ( 1- 
\chi ) v \,.
\end{equation}
If we now put 
\[ E  ( z) \stackrel{\rm{def}}{=}
 E_0 ( z) - E_+ ( z) R_+ ( z) E_0 ( z) \,, \ \ 
E_- ( z)  \stackrel{\rm{def}}{=} 
E_{ 0 , - } ( z) v - E_{-+ } ( z) R_+ ( z) E_0 ( z) \,,\]
then $ {\mathcal E} ( z ) $ given by \eqref{eq:E} is a right
microlocal inverse of $ {\mathcal P } ( z) $. 

To show that it is also a left inverse, we observe that 
\[ {\mathcal P} ( z ) ^* = 
\left( \begin{array}{ll} 
 - ( i / h ) P ( z) & R_+ ( z )^* \\
\; \ R_- ( z )^* & \ \ 0 \end{array} \right) \; : \; 
{\mathcal D}' (  X ) \times {\mathcal D}' ( \RR^{n-1}) 
\; \longrightarrow \;
{\mathcal D}' (  X ) \times {\mathcal D}' ( \RR^{n-1}) \,,\]
is microlocally defined in the same region as $ {\mathcal P} ( z) $
and is essentially of the same form but with $ W_+ $ replaced by $ W_-$
and $ \chi $ by $ 1 - \chi $:
\[ \begin{split}
& R_+ ( z) ^* = [( i/h) P( z) , \chi]_{ W_+}  K(z) \,, \\
& R_- ( z) ^* = K_b ( z) ^* [( i / h ) P ( z) , \chi] _{ W_-} \,.
\end{split}
\]
To see this we first note that 
\[ K_b ( z) ^* [( i / h ) P ( z) , \chi] _{ W_-} K_b ( z)  = - Id \,.\]
In fact, as in the proof of Proposition \ref{p:3.2}, \eqref{eq:4.2}
is invariant under the change of $ \chi $ and $ W_\pm $, as 
long as \eqref{eq:4.chi} hold. In particular, for
$ 0 < t < T ( 0 ) - \epsilon $, 
\[ K_b( z)^* [ ( i / h ) P ( z) , (\exp ( tH_{p(0)} )^* \chi]_{
\exp ( - t H_{ p( 0 )} ) W_+} K_b( z ) = Id \,. \]
For $ t \sim T(0)/2 $, $ W_+ $ is moved to $ W_- $, and $
 (\exp ( - tH_{p(0)} )^* \chi $ satisfies the properties of $ 1 - \chi $.
Hence, using the idependence of $ \chi$,
\[  K^b ( z) ^* [ ( i / h ) P ( z) , ( 1 - \chi)  ]_{W_- } K_b ( z) = 
Id_{ {\mathcal D}' ( \RR^{n-1} ) }
\,, \ \ \text{ microlocally near $ ( 0 , 0 ) \in T^* \RR^{n-1} $.} \]
If we now replace $ K ( z) $ by $ K_b ( z) $, then $ K( z) $ plays the 
r{\^o}le of $ K_b ( z) $, and this proves that $ R_+( z) ^* $ is the 
same as $ - R_-( z)$ with $ W_+ $ and $ W _- $ switched and $ \chi $
replaced by $ 1 -\chi $.

Hence, a similar argument to the one used for the construction of 
$ {\mathcal E} ( z) $ shows that $ {\mathcal P}( z) ^*$ has 
a right inverse, 
\[  {\mathcal F} ( z)^* = \left( \begin{array}{ll}
F (z ) & F_+ ( z ) \\
F_- ( z) & F_{ - + } ( z) 
\end{array} \right)^* \,. \]
Then $ {\mathcal F} ( z) $ is a left inverse of $ {\mathcal P } ( z) $,
and the usual argument ($ {\mathcal F} ( z) = {\mathcal F}( z) 
{\mathcal P}( z) {\mathcal E} ( z) = 
{\mathcal E }( z) $, microlocally) shows that it is equal to 
our right inverse. 

\medskip
\noindent
{\bf Remark.} By constructing part of the left inverse directly we 
can arrive at a simpler expression for $ E_- ( z)$:
\begin{equation}
\label{eq:4.em}
 E_- ( z ) = - (M ( z) K_f( z) ^* \chi +  K_b ( z ) ^* ( 1 - \chi )) \,,
\end{equation}
and it is useful to have it. To obtain it
we will directly solve the problem 
\begin{equation}
\label{eq:4.15} 
 \left\{ \begin{array}{l}   \widetilde E_- ( z ) ( i /h ) P ( z) 
+ E_ {-+  } ( z ) R_+ ( z ) = 0 \\
\widetilde E_- ( z) R _- ( z ) = Id_{ {\mathcal D}' ( \RR^{n-1} )} 
\end{array} \right. \end{equation}
Motivated by the structure of $ E_+ (z)$ and the fact that $ R_- ( z) $
is close to being an adjoint of $ R_+ ( z) $ (if it were, then 
$ E_- ( z ) $ would simply be the ajoint of $ E_+ ( z ) $), we put 
\begin{equation}
\label{eq:4.14}
\widetilde E_- ( z ) = - (M ( z) K_f( z) ^* \chi +  K_b ( z ) ^* ( 1 - \chi )) \,. 
\end{equation}
We now compute
\[ \begin{split} 
- \widetilde E_- ( z ) R_- ( z) &  = 
\left(M ( z) K_f( z) ^* \chi +  K_b ( z ) ^* ( 1 - \chi ) \right) 
 [ ( i / h )P ( z ) , \chi ]_{W_-} K_b ( z) \\
& = K^b ( z) ^* [ ( i / h ) P ( z) , \chi ]_{W_- } K_b ( z) \,.
\end{split}
\]
To analyze the last expression, we note that $ K ( z) $, in the 
definition of $ K_\bullet ( z ) $ was chosen, in Lemma \ref{l:4.1},
so that $ K( z ) ^* [ ( i / h ) P ( z) , \chi]_{W_+} K( z ) = Id $.
As in the proof of Proposition \ref{p:3.2}, 
this is invariant under the change of $ \chi $ and $ W_\pm $, as 
long as \eqref{eq:4.chi} hold: for 
$ 0 < t < T ( 0 ) - \epsilon $, 
\[ K_b( z)^* [ ( i / h ) P ( z) , (\exp ( tH_{p(0)} )^* \chi]_{
\exp ( - t H_{ p( 0 )} ) W_+} K_b( z ) = Id \,. \]
For $ t \sim T(0)/2 $, $ W_+ $ is moved to $ W_- $, and $
 (\exp ( - tH_{p(0)} )^* \chi $ satisfies the properties of $ 1 - \chi $.
Hence, using the independence of $ \chi$,
\[  K_b ( z) ^* [ ( i / h ) P ( z) , ( 1 - \chi)  ]_{W_- } K_b ( z) = 
Id_{ {\mathcal D}' ( \RR^{n-1} ) }
\,, \ \ \text{ microlocally near $ ( 0 , 0 ) \in T^* \RR^{n-1} $.} \]
This shows that $ \widetilde E_ - ( z ) R_- ( z ) = Id $ and we need to 
verify the first identity in \eqref{eq:4.15}. For that we use 
$ K_\bullet ( z) ^* P ( z) = 0 $, $ M ( z) K_f ( z ) ^* = K_b ( z ) ^* $,
near $ { ( 0 , 0) } \times W_-  \subset T^* \RR^{n-1} \times 
T^* X $, to obtain 
\[ \begin{split}
- \widetilde E_- ( z) ( i /h ) P ( z) & =  ( M ( z) K _f ( z)^* \chi + 
K_b ( z ) ^* ( 1 - \chi ))  ( i / h ) P ( z) \\
& = M ( z) K_f ( z )^* [ \chi , ( i / h ) P ( z) ]_{ W_+ } - 
K_b ( z)^* [\chi , ( i / h) P ( z) ]_{W_+} \\
& = K( z) ^* [ ( i / h ) P ( z) , \chi ] _{W_+} - 
M ( z) K( z) ^* [ ( i /h ) P ( z) , \chi ]_{ W_+} \\ 
& = ( 1 - M ( z) ) R_+ ( z ) = E_{-+}( z ) R_+ ( z ) \,, 
\end{split}
\]
and that establishes \eqref{eq:4.15}, so $ \widetilde E_- ( z) 
= E_- ( z) $ and we have \eqref{eq:4.em}.

\medskip

So far we considered only the case of $ z \in \RR $, and 
$ P ( z ) = P( z)^*$. Arguing as at the end of Sect.\ref{qt}, we see that
all the operators occuring in the construction of $ {\mathcal P }( z) $
and $ {\mathcal E} (z ) $ have almost analytic extensions to 
$ |\Im z | < L h \log ( 1/h )$ for any $ L$. It follows that the 
extention of $ {\mathcal E} (z )$ is a microlocal 
 inverse of the extension of $ {\mathcal P } ( z)$
modulo $ |\Im z |^\infty $, which in this neighbourhood of the real 
axis is $ {\mathcal O} ( h^\infty) $, that is, it remains a 
microlocal inverse. The bounds on the continuation of $ {\mathcal E} 
( z ) $ follow from \eqref{eq:2.4}.
This gives
\begin{prop}
\label{p:4.2}
Let $ P ( z) $ be an almost analytic extension of the self-adjoint
family of operators $ P ( z) \in \CI ( I_z ; \Psi^{ 0, k} ( X ) ) $,
such that 
\begin{gather*}
  \text{ The flow of $ H_p $ has a closed orbit $ \gamma $, }\\
\text{ on which $  p = \sigma ( P ( 0 )) = 0  $ and 
$  dp \neq 0$.}  \\
\end{gather*}
Then, there exist operators $ R_\pm ( z) $, defined in $ |\Im z | 
\leq L h \log ( 1/h) $, such that 
\[  {\mathcal P} ( z ) = \left( \begin{array}{ll} 
 ( i / h ) P ( z) & R_- ( z ) \\
\; \ R_+ ( z ) & \ \ 0 \end{array} \right) \; : \; 
{\mathcal D}' (  X ) \times {\mathcal D}' ( \RR^{n-1}) 
\; \longrightarrow \;
{\mathcal D}' (  X ) \times {\mathcal D}' ( \RR^{n-1}) \,,\]
defined microlocally near $ \gamma \times {(0, 0)} $, has a microlocal
inverse there:
\begin{gather*}
 {\mathcal E} ( z) = \left( \begin{array}{ll}
E (z ) & E_+ ( z ) \\
E_- ( z) & E_{ - + } ( z) 
\end{array} \right)  = 
{\mathcal O} ( e^{ C |\Im z | / h } ) : L^2 ( X ) \times L^2 ( \RR^{n-1} ) 
\rightarrow L^2 ( X ) \times L^2 ( \RR^{n-1} )  \,, \\
\bar \partial_z {\mathcal P } ( z) = {\mathcal O} ( |\Im z |^\infty )
\end{gather*}
and   $ E_{ -+} ( z) = I - M ( z) $, where $ M ( z) $ is the quantum 
monodromy operator defined by \eqref{eq:3.mon}. 
\end{prop}

\medskip
\noindent
{\bf Remark.} The constant $ C $ in the estimate of the norm of 
$ {\mathcal E} ( z) $ could be described more explicitely if stronger
conditions on $ P ( z) $ were made. If we assumed \eqref{eq:5.ell} 
then $ C $ could be related to $ C_p $ in \eqref{eq:5.cp}.

\section{Proof of the trace formula}
\label{pf}

We can now prove the main result of the paper. We strengthen our 
assumptions further here by demanding that $ P ( z) $ is a smooth
family of operators, self-adjoint for the real 
values of the parameter, and elliptic off the real axis.
\begin{thm}
\label{t:1}
Let $ P (z ) \in \CI ( I_z ; \Psi_h ^{0, k} ( X ) ) $, 
$ I = ( -a , a ) \subset \RR $, be a family of self-adjoint, 
principal type operators, such that $ \Sigma_z = \{ m \; : \; 
\sigma ( P ( z) ) = 0 \} \subset T^* X $ is compact.
We assume that
\begin{gather}
\label{eq:5.ell}
\begin{gathered}
\sigma ( \partial_ z P ( z) ) 
\leq - C <  0 \,, \ \text{ near $ \Sigma_z $.}\\
| \sigma ( P ( z) ) | \geq C | \xi |^k \,, \ \text{ for $ |\xi | \geq C $.}
\end{gathered}
\end{gather}
We also assume  that for $ z $ near $ 0$, the Hamilton vector field, 
$ H_{p(z)} $, $ p(z) = \sigma ( P ( z ) )$, has a simple closed
orbit $ \gamma(z) \subset \Sigma_0 $ with perid $ T ( z) $, 
and that $ \gamma ( z) $ has
a neighbourhood $ \Omega $ such that 
\begin{equation}
\label{eq:t2}
 m \in \Omega \ \text{and} \ \exp t H_{p(z)} ( m ) = m \,, \
p ( m ) = 0 \,, \ 0 < |t| \leq T (z) N 
+ \epsilon \,, \ z \in I\,, 
\ \Longrightarrow \ m \in \gamma(z) \,, \end{equation}
where $ T (z) $ is the period of $ \gamma( z) $, assumed to depend smoothly on
$ z $. 
Let $ A \in \Psi^{0,0}_h ( 
X ) $ be a microlocal cut-off to a sufficiently small 
neighbourhood of $ \gamma(0) $.

Then if $ P( z) $ is 
 an almost analytic extension of $ P( z) $, $ z \in \RR $,
$ \chi \in \CIc ( I ) $, $ \tilde \chi \in 
\CIc ( \CC) $, its almost analytic extension, 
$ f \in \CI ( \RR ) $, and $ {\rm{supp}}\; \hat f \subset 
( - N ( C_p - \epsilon)   + C ,  N ( C_p -  \epsilon)  - C )  
\setminus \{ 0 \} $, we have,
\begin{gather}
\begin{gathered}
\label{eq:t1}
\frac1\pi  {\mathrm{tr}} \;  \int f ( z / h ) 
 \bar \partial_z \left[ \tilde \chi ( z )\;   \partial_z P (z) \;  P ( z) ^{-1} \right] A 
{\mathcal L} ( d z ) 
 = \\
- \frac{1 }{ 2\pi i } 
\sum_{  -N -1  }^{N -1 }  {\mathrm {tr}} \; 
\int_{\mathbb R} f (z/h) M(z, h)^k \frac{d}{dz} M ( z , h )
 \chi ( z) 
dz + {\mathcal O} ( h^\infty ) \,, \end{gathered} \end{gather}
where $ M ( z, h ) $ is the quantum monodromy operator defined in 
\eqref{eq:3.mon} with $ K ( z) $ satisfying \eqref{eq:4.2}.
The constant $  C_p > 0 $, in the condition on $ \hat f $ depends on 
$ p (z) $ only and is given in \eqref{eq:5.cp}. 
\end{thm}

We observe that the left hand side of \eqref{eq:t1} is independent of
the choice of the almost analytic extension of $ \chi $: if $ \tilde \chi^
\sharp $ is another extension then, then $ \tilde \chi - 
\tilde \chi^\sharp = {\mathcal O} ( |\Im z |^\infty )$. In view of
Lemma \ref{l:5.1} below, 
\[( \tilde \chi ( z ) - \tilde \chi^\sharp ( z) ) 
\;   \partial_z P (z) \;  P ( z) ^{-1} \]
is smooth in $ z $, and $ {\mathcal O} ( |\Im z |^\infty )\,. $
By Green's formula and holomorphy of $ f $, the corresponding integral
vanishes. 

As described in Sect.\ref{out} Theorem \ref{t:0} is an immediate 
consequence of Theorem \ref{t:1}.

Before proceeding with a proof we remark that we can assume that
\[ P ( z ) \in \Psi^{0,0}_h ( X) \,, \]
since $ P ( z) $ can be multiplied by an $z$-independent elliptic
$ B \in \Psi^{0, -k } _h ( X ) $, without changing \eqref{eq:t1}.

We start with a 
lemma which justifies taking the traces in \eqref{eq:t1}:
\begin{lem}
\label{l:5.1}
Under the assumptions of Theorem \ref{t:1}, $ P ( z) ^{-1} $ exists in 
$ U \setminus I $, where $ U $ is a complex neighbourhood of $ J  \Subset I 
$, and
\[  \| P ( z )^{-1} \| \leq  C |\Im z |^{-1} \,, \ \ 0<  
|\Im z | \leq 1/C  \,. \]
\end{lem}
\begin{proof}
Let $ \psi = \psi^w ( x , h D ; z ) $ be a microlocal cut-off to a
a small neighbourhood of $ \Sigma_z $. Let us put $ v = P ( z) u $, so 
that (semi-classical) elliptic regularity gives 
\begin{equation}
\label{eq:5.1.1}
\| ( 1 - \psi ) u \| \leq C \| v \| + {\mathcal O} ( h^\infty ) 
\| u \| \,.
\end{equation}
For complex values of $ z $ we write 
\[ P ( z ) = P ( \Re z ) + \Im z \;  Q ( z ) \,, \]
where $ P ( \Re z ) $ is self-adjoint and $ \sigma (Q ( z ))  > 1/C > 0 $
near $ \Sigma_z $. This shows that 
\begin{equation}
\label{eq:5.1.2} 
\Im \langle P( z ) \psi u , \psi u \rangle = \Im z \; \Re \langle 
Q ( z ) \psi u , \psi u \rangle \geq 
\Im z \; \left( \| \psi u \|^2 / C - {\mathcal O}(h^\infty ) \| u \|^2 
\right) \,,\end{equation}
where we used the semi-classical G{\aa}rding inequality. 

We also write
\[ 
\begin{split}
& 
\Im  \langle P ( z) u , u \rangle - \Im \langle P ( z) \psi u , \psi u 
\rangle = \Im z  \left( \langle Q ( z ) u , u \rangle - 
\langle Q ( z ) \psi u , \psi u \rangle \right) = \\
& \Im z \; {\mathcal O } ( 1) \| ( 1 - \psi ) u \| \| u \| = 
\Im z \; {\mathcal O} ( 1) \left( \| v \| \| u \| + {\mathcal O}
( h^\infty ) \| u \|^2 \right) \,, \end{split} \]
where we used elliptic regularity \eqref{eq:5.1.1} in the last 
estimate. Then, applying \eqref{eq:5.1.2},
\[ 
\begin{split}
\| u \| \|v \| & \geq  \Im \langle P ( z ) \psi u , \psi u \rangle 
- \Im z \; {\mathcal O} ( 1) \left( \| v \| \| u \| 
+ {\mathcal O} ( h^\infty ) \|u\|^2 \right) \\
& \geq \Im z  \left( \| \psi u \|^2/ C  - {\mathcal O}(1) \|v\| \| u \| 
- {\mathcal O} ( h^\infty ) \| u \|^2 \right) \,. \end{split} \]
For small $ \Im z $ the term $ \|v\| \| u \| $ on the left hand side 
can be absorbed in the right hand side, and by adding 
$ \Im z \| (1 - \psi ) u \|^2 $ to both sides we obtain
\[ \Im z \| u \|^2 / C  \leq  \| u \| \| v \| + {\mathcal O} ( h ^
\infty ) \Im z \| u \|^2  \, , \]
and that gives
\[ \| u \| \leq \frac{C}{ \Im z }  \| v \|   \,, \] 
proving the estimate for $ P ( z ) ^{-1} $.
\end{proof}

\medskip
\noindent
{\em Proof of Theorem:}
Using Proposition \ref{p:4.2} we can {\em formally} write
\[ P ( z) ^{-1} A = 
E ( z) A - E_+ (z ) E_{-+} ( z) ^{-1} E_- ( z) A \,, \ \ 
E_{-+ } ( z) = I  - M ( z) 
 \,, \]
microlocally near $ \Omega $, and for $ 0< |\Im z | \leq L h \log ( 1/h)$,
with any $ L $. 
To apply this formal expression rigourously,
we rewrite the left hand side of \eqref{eq:t1} as
\begin{gather}
\label{eq:5.1}\begin{gathered}
\frac{1}{\pi }  {\mathrm{tr}} \;  \int f ( z / h ) 
\bar \partial_z \tilde [ \chi ( z )\;  \partial_z P (z) \;  P ( z) ^{-1} ] \;
 A \;
{\mathcal L} ( d z ) 
= \\
\frac{1}{\pi}
\sum_{\pm}  {\mathrm{tr}} \;  \int_{\CC_{\pm}} f ( z / h ) 
\bar \partial_z \tilde [ \chi ( z )\;   \partial_z P (z) \;  P ( z) ^{-1} \;
]  A \;
{\mathcal L} ( d z ) \,. 
\end{gathered}
\end{gather}
Then, 
motivated by the formal Neumann series expansion of $ ( I - M( z ) )^{-1} $
we define
\begin{equation}
\label{eq:5.2}
T_N^+ ( z) \stackrel{\rm{def}}{=}
 E( z) A - E_+ ( z ) \sum_{ k=0}^N M(z)^k E_- ( z) A \,, 
\end{equation}
so that 
\begin{equation}
\label{eq:5.2'}
P(z)^{-1} A = T_N^+ ( z ) + P ( z) ^{-1} R_- ( z ) M( z)^{N+1} E_ - ( z) A
\,, \end{equation}
microlocally near $ \gamma $ and for $ 0 < |\Im z | \leq L h \log ( 1/h ) $,
for any $ L$.
In fact, from $ {\mathcal P} ( z ) {\mathcal E} ( z) = Id $, and
$ E_\pm ( z) = I - M ( z) $,  we have 
\[  P ( z) E_ + ( z ) = - R_ - ( z ) ( I -  M ( z) ) \,, \ \ 
P ( z) E( z) = I - R_- ( z) E_ - ( z) \,, \]
and hence
\[ \begin{split}
P ( z) T_N^+ ( z) & = P ( z) ( E ( z) A - E_+ ( z) 
\sum _{ k=0}^N M( z) ^k E_- ( z) A )  \\
& = A - R_- ( z) E_- ( z) A + R_-( z) ( I - M ( z) ) \sum_{ k = 0 }^N 
M( z) ^k E_ - ( z) A  \\
& = A - R_- ( z) M( z)^{N+1} E_- ( z) A \,, \end{split} \]
which gives \eqref{eq:5.2'}.

To use this in \eqref{eq:5.1} we need to have the support of the
almost analytic extension of the cut-off function $ \chi $ to be
contained in the region where $ |\Im z | \leq L h \log (1/h ) $. 
To do that we follow the method of \cite[Sect.12]{DiSj} by fixing
an almost analytic extension of $ \chi$, $ \chi^\# $, and 
then putting
\[ \tilde \chi = \tilde \chi_{L,h} = \chi^\# \psi_{L,h} \,, \ \ 
\psi_{ L, h } (z ) = \psi \left( \frac{\Im z}{ L h \log( 1/h)} \right)\,,
\ \ \psi( t) = \left\{ \begin{array}{ll} 
1 \,, & |t | < 1/2 \\
0  \,, &  |t | > 1 \end{array} \right.
\,. \]
By the remark after the statement of Theorem \ref{t:1}, \eqref{eq:t1} is
independent of the choice of $ \tilde \chi $ and hence we can use 
$ \tilde \chi _{ L , h } $. 
Since now $ {\mathcal O} ( |\Im z |^\infty ) = {\mathcal O } ( h ^ \infty ) 
$, the almost analyticity of $ 
P ( z) $ also shows that the left hand side of \eqref{eq:t1} can
be rewritten as
\begin{equation}
\label{eq:t1'}
\frac1\pi  {\mathrm{tr}} \;  \int f ( z / h ) 
 \bar \partial_z  \tilde \chi ( z )\;   \partial_z P (z) \;  P ( z) ^{-1} 
\;  A 
{\mathcal L} ( d z ) \,,
\end{equation}
and this is what we will use from now on.

We claim that with the choice of $ \tilde \chi $ above
\begin{equation}
\label{eq:5.3}
\begin{split}
& {\mathrm{tr}} \;  \int_{\CC_{+}} f ( z / h ) 
\bar \partial_z \tilde \chi ( z )\;  \partial_z P (z) \;  P ( z) ^{-1} \; 
A \;
{\mathcal L} ( d z ) = \\
& \ \ 
{\mathrm{tr}} \;  \int_{\CC_{+}} f ( z / h ) 
\bar \partial_z \tilde \chi ( z )\;   \partial_z P (z) \;
 T_N^+( z ) {\mathcal L} ( d z ) + {\mathcal O} ( h^{ L/ C} ) \,, 
\end{split}
\end{equation}
where $ C $ is fixed depending on $ N $ and $ \supp \hat f $.

To show this we first need the following
\begin{lem}
\label{l:5.2}
The almost analytic continuation of the monodromy operator satisfies, 
for $ z $ sufficiently close to $ 0 $, and
for
any $ L $, 
\begin{equation} \begin{split}
& \| M ( z ) \| \leq e^{ - ( C_p - \epsilon ) \Im z /h } + {\mathcal 
O} (h^\infty ) \, , \ \ 
0 < \Im z < L h \log ( 1/ h )  \,, 
\\
& \| M ( z )^{-1} \| \leq e^{ ( C_p - \epsilon ) \Im z /h }
+ {\mathcal O} (h^\infty )  \, , \ \ 
- L h \log ( 1/h ) < \Im z <  0 \,, 
\\ &  \ \ C_p = - \int_0^{T(0)} \sigma ( \partial_z P ( z) ) ( \exp
( t H_{p( 0)} ) ( m_0 ) dt \,,  \label{eq:5.cp}
\end{split}
\end{equation}
where $ \epsilon > 0 $ can be taken arbitrarily small by shrinking the
neighbourhood of $ \gamma $.
The constant $ C_p $ is positive thanks to \eqref{eq:5.ell}.
\end{lem}
\begin{proof}
We use the differential equation \eqref{eq:3.m} and observe that for
$ z $ real, and
$  m_0 ( z) \in \gamma ( z )$, 
\[ \sigma ( K ( z) ^{-1} ( \QC ( z) - Q ( z ) ) K  ( z) ) ( 0 , 0 ) 
= - \int_0^{T(z)} \sigma ( \partial_z P ( z) ) ( \exp
( t H_{p( z)} ) ( m_0(z)  ) dt \,.\]
Hence, writing $ z = x + i y $, $ 0 < y < L h \log( 1/h ) $,
and $ B ( z ) = K ( z) ^{-1} ( Q ( z) - \QC ( z ) ) K  ( z) $,  we have,
for $ v \in {\mathcal D}' ( \RR^{n-1} )$, with $ WF ( v ) $ close to 
$ ( 0 , 0 )$, 
\[ \begin{split}
h \frac{d}{dy} \left( \| M ( z) v \|^2 \right) 
& = h \frac{d}{dy} \langle M ( z ) v , M( z  ) v \rangle \\
& = i h ( \partial_z - \partial_{\bar z })
 \langle M ( z ) v , M( z  ) v \rangle \\
& = 
- \langle B ( z) M ( z) v , M( z)v \rangle - \langle M ( z ) v , B( z) 
M ( z) v \rangle + {\mathcal O} ( |\Im z |^\infty ) \| v \|^2 \\
& = - \langle ( B ( z) + B( z)^* ) M ( z ) v , M ( z) v \rangle + 
{\mathcal O} ( |\Im z |^\infty ) \| v \|^2 \,.
\end{split}
\]
The  G{\aa}rding inequality now shows that for $ x $ small 
enough,
\[ h \frac{d}{dy} \left( \| M ( z) v \|^2 \right)  \leq 
- ( C_p - \epsilon ) \| M ( z) v \|^2 + {\mathcal O} (y ^\infty ) \|v\|^2 
\,. \]
Since by Proposition \ref{p:3.2}, $ \| M ( x ) v \|^2 = \| v\|^2 ( 1+ 
{\mathcal O}( h^\infty )) $, the lemma follows.
\end{proof}

\medskip
\noindent
{\em Proof of \eqref{eq:5.3}:}
By \eqref{eq:5.1} and \eqref{eq:5.2'}
we need to estimate 
\begin{equation}
\label{eq:5.3'}  {\mathrm{tr}} \;  \int_{\CC_{+}} f ( z / h ) 
\bar \partial_z \tilde \chi_{L, h}
 ( z )\;   \partial_z P (z) \; P ( z) ^{-1} \;
 R_- ( z) M( z) ^{N+1} E_ -  ( z) \; A \;
{\mathcal L} ( d z )  \,, \end{equation}
where by Lemmas \ref{l:5.1} and \ref{l:5.2} we have
\[ \begin{split}
& \| M( z)^{N+1}\| \leq  \| M ( z) \|^{N+1} \leq 
e ^ { - ( C_p - \epsilon ) ( N + 1) \frac{ \Im z}{h}} \,, \ \ 
0 \leq \Im z \leq L h \log (1/h) \,, \\
& \| P ( z) ^{-1} \| \leq 1 / |\Im z| \,.\end{split}
\]
All the operators coming from $ {\mathcal P }( z) $ and 
$ {\mathcal E}( z ) $ are bounded by $ \exp ( C | \Im z | / h  )$, 
and if $ \supp \hat f \subset 
[ -b + C , b - C ] $, then 
\[ |  f( z     / h ) | \leq C e^{ ( b - C ) \frac{ |\Im z | }{h}} \,.\]
Using the definition of $\tilde \chi_{L,h} $, the above estimates, and
the characteristic function
\[ \rho_{ L, h } ( t) =  
\bbbone_{ L h \log ( 1/h ) /2 \leq t \leq Lh \log ( 1/h ) } \,,\]
we can bound \eqref{eq:5.3'} by a constant times
\[  h^{-n} \int
\left( |\bar \partial_z \chi^\# ( z ) | + 
 ( L h \log ( 1/h ) )^{-1} \rho_{L, h} ( \Im z ) \right) 
|\Im z |^{-1} e^{ \frac{\Im z}{h} ( b - (C_p - \epsilon) ( N + 1) )} 
{\mathcal L} ( dz ) \leq \]
\[ C h^{-n}\int_{ 0 \leq \Im z < L h \log( 1/h) } 
  |\Im z |^M e^{ \frac{\Im z}{h} ( b - (C_p - \epsilon) ( N + 1) )} 
{\mathcal L} ( dz )  \]
\[ + \; 
C h^{-n} ( L h \log ( 1/h ) )^{-2} h^{ (( C_p - \epsilon) ( N + 1 ) - b)L }
\leq  C'_N h^{ L/C -n } \,, 
\]
where $ C> 0 $ is fixed. 
\stopthm

With \eqref{eq:5.3} established, we have to 
study the leading term on its right hand side 
which we rewrite using the definition \eqref{eq:5.2} and
the cyclicity of the trace:
\begin{gather}
\label{eq:5.ga}
\begin{gathered}
\frac{1}{\pi} {\rm{tr}}\; \int_{\CC_+} 
\bar \partial_z ( \tilde \chi_{L,h} ) \partial_z P ( z) \; E( z)\;  A \;
f ( z/h )  {\mathcal L} ( dz)
- \\
\frac{1}{\pi} {\rm{tr}}\; \int_{\CC_+} 
\bar \partial_z ( \tilde \chi_{L,h} ) \;
\left( \sum_{ k = 0}^N M ( z) ^k E_ - ( z) A \partial_z P ( z) 
 E_+ ( z) \right) 
f ( z / h ) {\mathcal L} ( dz) \,. 
\end{gathered}
\end{gather}
Since all the operators are almost analytic 
(in particular $ 
\bar \partial_z E_{\bullet}( z) = {\mathcal O} ( h^\infty ) $
on the support of $ \tilde \chi _{ L , h } $)
and $ f ( z/ h )$ is holomorphic, we can apply Green's formula
and reduce the first integral  to an integral over the real axis:
\[ \frac{1}{2 \pi i }  {\rm{tr}}\;  \int_{\RR} \chi ( z ) f ( z / h  ) 
\; \partial_z P ( z) E( z) A d z  
\,. \]
To analyze the second term (with integration still over $ \CC_+ $) we see
that the explicit expressions $ E_- ( z) $ and $ E_+ ( z ) $, \eqref{eq:4.7}
and \eqref{eq:4.em}, show that
\[  E_- ( z ) A \partial_z P ( z) E_+ ( z ) = A_1 ( z ) M ( z ) + A_2 ( z ) 
\,, \ \ A_j ( z ) \in \CI ( I_z ; \Psi
^{0, - \infty } ( \RR^{n-1} ) )\,. \]
To analyze the contributions to the trace we need the following simple
\begin{lem}
\label{l:4.3}
Suppose that $ A \in \Psi_h^{ 0 , - \infty } ( \RR^{n-1} ) $ and
$ U \in I^0 _h ( \RR^{n-1} \times \RR^{n-1}; { C}' )$, 
satisfy
\[  m \in WF_h ( A ) \ \Longrightarrow ( m ,  m ) \notin { C} \,.\]
Then, for $ \chi \in \CIc ( \RR^{n-1} )$,
\[  \tr \chi A U  = {\mathcal O} ( h^{\infty } ) \,. \]
\end{lem}
\begin{proof}
It is clear that if $ \chi_1, \chi_2 \in \CIc (T^* \RR^{n-1} ) $ have 
disjoint supports then $ \tr \chi_1 A U \chi_2 = {\mathcal O} ( h^\infty ) $,
where we denoted the quantizations by the same symbols. Using the 
hypothesis, we can write $ \chi A U $  as a sum of 
negligible terms ($ {\mathcal O} ( h^\infty ) $), and of terms of that
form.
\end{proof}

The assumption \eqref{eq:t2} implies that 
the main contribution (modulo $ {\mathcal O} ( h^\infty ) $
as usual) to the trace of 
\[ M ( z) ^k E_- ( z ) A \partial_z P ( z ) E_+ ( z ) \,, \ \ k > 0 \,, \]
comes from an arbitrarily small neighbourhood of the fixed point, $ ( 0 , 0 )
\in T^* \RR^{n-1} $, 
of $ C (z)$, the canonical relation of $ M ( z ) $. We can therefore 
replace $ A $ by $ 1 $ and introduce a microlocal cut-off, $
\rho^w $,  to a neighbourhood
of $ ( 0 , 0 ) $:
\[ \tr  M ( z) ^k E_- ( z ) A \partial P ( z ) E_+ ( z ) 
=  M ( z) ^k E_- ( z )  \partial P ( z ) E_+ ( z ) \rho^w + {\mathcal O} 
( h^\infty )   
\,, \ \ k > 0 \,. \]
For $ k = 0 $ the same discussion is valid for the contribution of 
$ M ( z) A_1 ( z ) $ in $  E_- ( z )  \partial P ( z ) E_+ ( z )     $, but 
for the pseudo-differential contribution, $ A_2 ( z ) $, we need to use
the support assumption on $ \hat f $: $ 0 \notin \supp \hat f $. We write
\[
\tr \int_{ \CC_+ } f ( z / h ) \bar \partial_z \tilde \chi ( z) 
A_2 ( z )  \rho^w {\mathcal 
L } ( dz ) = 
\tr \int_{ \RR } f ( z / h ) \chi ( z) A_2 ( z )  \rho^w  dz  = 
{\mathcal O} ( h^\infty ) \,, \]
by the standard argument: put $ g ( z ) = \tr A_2 ( z) \rho^w $, 
so that by Plancherel's theorem
\begin{equation}
\label{eq:plan} \begin{split}
2 \pi \int_{ \RR } f \left (\frac{z}{h} \right) g ( z) dz & = 
h \int_{\RR }  \hat f ( h \zeta ) \hat g ( - \zeta ) d\zeta\\
& = 
{h^{M+1}}  \int_{\RR }   \hat f ( h \zeta ) / ( h \zeta)^M
\zeta^M  \hat g ( - \zeta ) d\zeta \\
& = {\mathcal O} ( h^{M+1}  ) \,. \end{split} \end{equation}

Hence the second term in \eqref{eq:5.ga} becomes 
\begin{equation}
\label{eq:3sj}
- \frac{1}{\pi} {\rm{tr}}\; \int_{\CC_+} 
\bar \partial_z  \tilde \chi \;
\left( \sum_{ k = 0}^N M ( z) ^k E_ - ( z)  \partial_z P ( z) 
 E_+ ( z) \right) \rho^w  
f ( z / h ) {\mathcal L} ( dz) \,, \end{equation}
where $ \rho^w $ is a microlocal cut-off to a neighbourhood of $ ( 0 , 0 )$.

We recall that when $ R_\pm $ are
independent of $ z $, the following standard formula holds:
\[ \partial_z M  ( z) = -  \partial_z E_{ -+} ( z)  = 
E_- ( z) \; \partial_z P ( z ) \; E_+ ( z) 
\,, \]
as is easily seen from $ \partial_z {\mathcal E} = 
- {\mathcal E} \partial_z {\mathcal P } {\mathcal E} $.
In the general case, the same argument gives
\begin{equation}
\label{eq:5.dm}
\begin{split}
 \partial_z M  ( z) & = -  \partial_z E_{ -+} ( z) \\
& = 
 E_- ( z) \; \partial_z P ( z ) \; E_+ ( z) + E_{ - + } ( z) 
\; \partial_z R_+ ( z ) \; E_+ ( z) + E_- ( z ) \;
\partial_z R_- ( z ) \;  E_{ - + } ( z ) \,. \end{split}
\end{equation}

Inserting this we obtain the following expression for \eqref{eq:3sj}:
\begin{equation}
\label{eq:4sj}
\begin{split}
& - \frac{1}{\pi} \tr \int_{\CC_+} f ( z / h ) \bar \partial \tilde 
\chi ( z ) \sum_{k=0}^N M ( z) ^ k \partial_z M ( z) \rho^w {\mathcal L}
( dz ) \\
& + \frac{1}{\pi} \tr \int_{\CC_+} f ( z / h ) \bar \partial \tilde 
\chi ( z ) \sum_{k=0}^N M ( z) ^ k ( 1 - M ( z))^{-1}
 \partial_z R_+ ( z) E_+ ( z) 
\rho^w {\mathcal L}
( dz ) \\
& + \frac{1}{\pi} \tr \int_{\CC_+} f ( z / h ) \bar \partial \tilde 
\chi ( z ) \sum_{k=0}^N M ( z) ^ k E_- ( z) \partial_z R_- ( z) 
( 1 - M ( z))^{-1}  \rho^w M ( z) ^ k 
 {\mathcal L}
( dz ) \\
& = J_1 + J_2 + J_3 \end{split}
\end{equation}
By Green's formula
\[ J_1 = - \frac{1}{ 2 \pi i } \sum_{k=0}^N \tr \int_{\RR} 
f(  z / h )  \chi ( z) M ( z ) ^ k \partial_z M ( z ) \rho^w dx \,, \]
which is a term appearing in \eqref{eq:t1}. We want to show that
the remaining two terms, $ J_2$, $J_3 $, are negligible.

To see this we need
\begin{lem}
\label{l:5.f}
We have
\[ \partial_z R_+ ( z) E_+ ( z) \,, \ \ E_- ( z ) \partial_z R_- ( z) 
\in \CI ( I_z ; \Psi^{ 1, - \infty }_h  ( X ))  \,. \]
\end{lem}
\begin{proof}
From the definitions \eqref{eq:4.3}, \eqref{eq:4.7}, and from \eqref{eq:4.2}
we see that 
\[ \partial_z R_+ ( z ) E _+  ( z ) = 
\partial_z \left( K^* ( z ) [ ( i / h ) P ( z ) , \chi] _{W_+}  \right) 
K ( z) = - 
 K^* ( z ) [ ( i / h ) P ( z ) , \chi] _{W_+} \partial_z K ( z) \,.\]
From the proof of Proposition \ref{p:3.m} we recall that
\[ h D_z K ( z ) = - Q ( z) K ( z) \,, \]
and hence
\[ \partial_z R_+ ( z ) E _+  ( z ) =
( i / h )  K^* ( z ) [ ( i / h ) P ( z ) , \chi] _{W_+} Q ( z) K ( z) \,.\]
This expression is microlocal near $ ( 0 , 0 ) $ and as far as $ K ( z) $ 
is concerned microlocal near $ ( m_0 , ( 0 , 0 ) ) $. Hence we can use
a model given in Proposition \ref{p:2.1}: $ P ( z ) = h D_{x_1} $ 
(the microlocal $z$-dependent conjugation will not affect the 
uniform pseudo-differential behaviour), and
\[  K ( z) u ( x_1, x' ) = \frac{1}{( 2 \pi h )^{n-1} } 
\int e^{ i ( \langle y' , \eta' \rangle  - \phi_z ( x' , \eta' ) )/ h }
a_z ( x' , \eta' ) u ( y') dy' d\eta' \,, \]
where we used local representation of the $h$-Fourier Integral Operators
(see the proof of Proposition \ref{p:6.2} below for the derivation of a local
representation). After composing the operators, and applying the stationary
phase method we arrive at the following expression for the kernel of 
$ \partial_z R_+ ( z ) E _+  ( z )$:
\[ \frac{1}{( 2 \pi h )^{n-1} } 
\int e^{ i (  \phi_z ( x' , \eta' ) - \phi_z ( y', \eta' )) / h  }
A_z (y',  x' , \eta' ) d\eta' \,, \ \ A_z \in S^{ 1, - \infty }  \,,\]
which by a standard ``Kuranishi trick'' argument (see the appendix)
shows that we get
a smooth $ z $-dependent family of pseudo-differential operators.
\end{proof}

In $ J_2 $ we can replace $   \sum_{k=0}^N M ( z) ^ k ( 1 - M ( z))$ 
by $ 1 - M ( z )^{N+1} $. As in the proof of \eqref{eq:5.3}, we show
that the term corresponding to $ M ( z)^{N+1} $ is negligible. The
remaining term is transformed to an integral over $ \RR $:
\[  \frac{1}{ 2 \pi i } \int_{\RR} 
f(  z / h )  \chi ( z) \tr ( \partial_z R_+ ( z) E_+ ( z)  \rho^w ) dx \,,
\ \ 0 \notin \supp \hat f 
\]
which is negligible by Lemma \ref{l:5.f} and \eqref{eq:plan}. Similar arguments
then apply to $ J_3 $. 

To summarize, we have shown that 
\begin{gather*}
 {\mathrm{tr}} \;  \int_{\CC_{+}} f ( z / h ) 
\bar \partial_z \tilde \chi ( z )\;  \partial_z P (z) \;  P ( z) ^{-1} \; 
A 
{\mathcal L} ( d z ) = \\
 \frac{1}{2 \pi i }  {\rm{tr}}\;  \int_{\RR} \chi ( z ) f ( z / h  ) 
\; \partial_z P ( z) E( z) A d z  
- \frac{1}{ 2 \pi i } \sum_{k=0}^N \tr \int_{\RR} 
f(  z / h )  \chi ( z) M ( z ) ^ k \partial_z M ( z ) \rho^w dx + 
{\mathcal O} ( h^\infty ) \,, 
\end{gather*}

We proceed in a similar way for the integral over $ \CC_- $ in \eqref{eq:5.1}.
We write $  I - M ( z)  =  - M ( z) ( I - M ( z)^{-1} ) $, and motivated by
the resulting formal Neumann series put
\[ T_N^- ( z) \stackrel{\rm{def}}{=}
 E( z) A +  E_+ ( z ) \sum_{ k=1}^N ( M(z)^{-1} ) ^k E_- ( z) A \,.\]
The same arguments apply and Green's formula gives 
\[ - \frac{1}{2 \pi i }  {\rm{tr}}\;  \int_{\RR} \chi ( z ) f ( z / h  ) 
E( z) A d z  -  
\frac{1}{2 \pi i }  {\rm{tr}}\;  \int_{\RR}  
\chi ( z) f (z / h ) 
\sum_{ k = 1}^{N +1 } M ( z) ^{-k}  \partial_z M ( z ) 
 dz + {\mathcal O} ( h^\infty ) \,. \]

When we now add the contributions from the integrations over $ \CC_\pm $
we see that the integrals involving $ E ( z ) $ cancel and the remaning
terms give \eqref{eq:t1}

\stopthm

\section{Trace formula for non-degenerate closed trajectories}
\label{tr}

We say that a  closed trajectory $ \gamma ( z )$ of $ P ( z ) $
{\em $ N$-fold non-degenerate} if 
\begin{equation}
\label{eq:6.non}
 \det ( I - (d C ( z ) _{m_0(z)})^k  ) \neq 0 \,, \ \ 
0 \neq |k | \leq N \,, \end{equation}
where $ C ( z ) $ is the Poincar{\'e} map of $ \gamma ( z )$, 
\eqref{eq:4.poin}.
When this holds our theorem translates into the standard 
semi-classical trace formula, generalizing the 
Gutzwiller, Balian-Bloch, and Duistermaat-Guillemin trace formul{\ae}.

We start by a general discussion of traces of Fourier 
Integral Operators.

\begin{lem}
\label{p:6.1}
Suppose that $ B $, microlocally defined near $ ( 0 , 0 ) \in 
T^* \RR^n $ is given by 
\begin{equation}
\label{eq:6.1}
Bu(x)={1\over (2\pi h)^n} \iint e^{i(\phi (x,\eta )-y\eta
)/h}b(x,\eta ;h)u(y)dyd\eta 
\end{equation}
where $\phi (x,\eta )$ is defined near $ (0,0)$,  $a$ is
a classical symbol of order $ 0 $, supported near
$ ( 0 , 0 ) $, $\phi '(0,0)=0$, and 
 $ \phi''_{x \eta } ( 0 , 0 ) \neq 0 $. 
The corresponding canonical transformation is given by 
\begin{equation}
\label{eq:6.2}
{\kappa :(\phi '_\eta(x,\eta ) ,\eta )  \mapsto (x,\phi
'_x(x,\eta )) \,. }\end{equation}
It is defined 
between two neighbourhoods of $ (0,0)$  and we assume that
$ (0,0) $ is its only fixed point there, and
\begin{equation}
\label{eq:4}
\det (d\kappa (0,0)- 1 ) \neq 0 \,. 
\end{equation}

Under these assumptions 
\begin{equation}
\label{eq:9}
 \tr B=i^{\frac 12 s }  {(b_0(0,0)+{\mathcal O}(h))e^{i\phi
(0,0)/h}\over
|{\det \phi ''_{\eta x}} \det (d\kappa (0,0)-1) |^{\frac12	} }\,, 
\ \ s = {\rm{sgn}} \; \left( \begin{array}{ll} \phi ''_{xx} &\phi
''_{x\eta }-1 \\ \phi ''_{\eta x}-1&\phi ''_{\eta \eta
}   \end{array} \right) \,,
\end{equation}
where the signature of a symmetric matrix $ A$, $ {\rm{sgn}}\; A $, is 
the difference between the number of positive and negative eigenvalues. 
\end{lem}
\begin{proof}
The fact that \eqref{eq:6.2} defines a smooth map is
equivalent to the assumption that
\begin{equation}
\label{eq:3}
\det \phi ''_{x \eta }\ne 0 \,. \end{equation}
Here and in the following, second derivatives of $\phi $
are computed at (0,0) if nothing else is specified.
The differential,
$d\kappa (0,0)$, is the map
$(\delta _y,\delta _\eta )\mapsto (\delta _x,\delta _\xi
)$, where
\[ \begin{split} & \delta _y=\phi ''_{\eta x}\delta _x+\phi
''_{\eta \eta} \delta _\eta \\
&  \delta _\xi =\phi
''_{x\eta }\delta _\eta +\phi ''_{xx}\delta _x \,. 
\end{split} \] 
Here we can express $\delta _x$ and $\delta _\xi $ in
terms of $\delta _y,\delta _\eta $ and it follows that
$d  \kappa (0,0)$ is given by the matrix:
\begin{equation}
\label{eq:5}
d \kappa (0,0)= 
\left( \begin{array}{ll} (\phi ''_{\eta x})^{-1}
& -(\phi _{\eta x}'') ^{-1} \phi _{\eta \eta }'' \\
\phi ''_{xx}(\phi ''_{\eta x})^{-1} & \phi ''_{x\eta }-
\phi _{xx}''(\phi ''_{\eta x})^{-1}\phi''_{\eta \eta }  
\end{array} \right)  \,.   \end{equation}
We find the following factorization:
\begin{equation}
\label{eq:6}
\begin{split}
& d\kappa (0,0)-1= \\
& \left( \begin{array}{ll} 
-(\phi _{\eta x}'')^{-1} & 0 \\ 0
& 1 \end{array} \right) 
\left( \begin{array}{ll} 0 &1 \\  1 &0 
\end{array} \right) 
\left( \begin{array}{ll} 1 & -\phi ''_{xx}(\phi
''_{\eta x})^{-1}\\ 0 &1 \end{array} \right) 
\left( \begin{array}{ll} \phi ''_{xx} &\phi
''_{x\eta }-1 \\ \phi ''_{\eta x}-1&\phi ''_{\eta \eta
}  \,. \end{array} \right) \,. 
\end{split}  \end{equation}
In particular,
\begin{equation}
\label{eq:7}
\det (d\kappa (0,0)-1)={1\over \det \phi ''_{\eta
x}}\det\left( \begin{array}{ll} 
\phi ''_{xx} &\phi ''_{x\eta }-1\\ \phi
''_{\eta x}-1 &\phi ''_{\eta \eta } \end{array} \right) \,. 
\end{equation}

Here 
$$ \left( \begin{array}{ll} \phi ''_{xx}&\phi ''_{x\eta }-1\\
\phi ''_{\eta x}-1&\phi ''_{\eta \eta } \end{array} \right) $$ 
is the Hessian
of $\phi (x,\eta )-x\eta $. The stationary phase method applied to 
the trace of \eqref{eq:6.1} gives
\[ 
{\rm tr\,}B =\left( \left(\det {1\over i} \left( \begin{array}{ll}\phi
''_{xx}&\phi ''_{x\eta }-1\\ \phi ''_{\eta x}-1&\phi
''_{\eta \eta } \end{array} \right) 
\right)^{-1/2} b_0 (0,0)+{\mathcal O}(h) \right)e^{i\phi
(0,0)/h}. \]
Here we choose the branch of the square root of the
determinant on the set of non-degenerate symmetric matrices
with non-negative real part 
which is equal to 1 for the identity. Using \eqref{eq:7}, we get
\eqref{eq:3}.
\end{proof}
 
To give a geometric meaning to the signature $ s $ appearing in 
\eqref{eq:9} in terms of a {\em Maslov index} we first
recall the definition of the H\"ormander-Kashiwara
index of a Lagrangian triple: let $ \lambda_1, \lambda_2, \lambda_3
$ be Lagrangian planes in a symplectic vector space $ ( V , \omega ) $,
and put
\begin{equation}
\label{eq:6.hk}
s ( \lambda_1 , \lambda_2 , \lambda_3 ) = {\rm sgn} \;
Q ( \lambda_1 , \lambda_2 , \lambda_3) \,, \end{equation}
where $ Q ( \lambda_1 , \lambda_2 , \lambda_3) $ is a quadratic form on 
$ \lambda_1 \oplus \lambda_2 \oplus \lambda_3 $ given by 
\[ 
Q ( \lambda_1 , \lambda_2 , \lambda_3) ( 
v_1 \oplus v_2 \oplus v_3 ) = \omega ( v_1, v_2 ) + \omega ( v_2, v_3) 
+ \omega ( v_3 , v_1 ) \,, \]
see 
\cite{LV}  for a comprehensive introduction. 
Here we only mention that if $ \lambda_i$'s are
mutually transversal, then $ s ( \lambda_1, \lambda_2 , \lambda_3 ) $
is the only symplectic invariant of such three Lagrangian planes. It is 
antisymmetric and satisfies the cocycle condition.

We then have 
\begin{lem}
\label{l:6.1}
Let $ V = T^* \RR^n \times T^*\RR^n $ with the symplectic form 
$ \omega = \omega_1 - \omega_2 $, where $ \omega_1 $ and $ \omega_2 $
are the canonical forms on the factors. In the notation of 
Lemma  \ref{p:6.1}, let $ \Gamma_{d \kappa } $ be the graph of 
$ d \kappa ( 0 , 0) $, $ \Delta \subset  T^* \RR^n \times T^*\RR^n $
be the diagonal, and $ M =  \{ 0 \} 
\oplus \RR^n  \oplus \RR^n \oplus \{ 0 \} 
\subset  T^* \RR^n \times T^*\RR^n $. Then
\[ s ( \Gamma_{d \kappa } , \Delta, M ) = - 
 {\rm{sgn}} \; \left( \begin{array}{ll} \phi ''_{xx} &\phi
''_{x\eta }-1 \\ \phi ''_{\eta x}-1&\phi ''_{\eta \eta
}   \end{array} \right) \,. \]
\end{lem}
\begin{proof}
Let us write
\[ \phi'' ( 0 , 0 ) = \left( \begin{array}{ll} \alpha & \beta 
\\ \beta^t & \gamma \end{array} \right) \,. \]
Then $ \Gamma_{ d \kappa } = \{ ( x , \alpha x + \beta \eta ;
 \beta^t x + \gamma \eta,  
\eta ) \; : \; ( x , \eta ) \in 
\RR^n \times \RR^n \} $. Since $ \Delta $ and $ M $ are
transversal, \cite[Lemma 1.5.4]{LV} says that
\[ s ( \Gamma_{ d \kappa } , \Delta , M ) = 
- {\rm{sgn}} \; \omega ( \pi \bullet , \bullet ) |_{\Gamma_{ d \kappa }}
\,, \]
where $ \pi : T^* \RR^n \times T^*\RR^n  \rightarrow M $ is the 
projection along $ \Delta $:
$ \pi ( x , \xi ;   y , \eta )  = ( 0 , \xi - \eta;  y - x , 0  ) $. 
In the $ ( x , \eta  ) $ coordinates on $  \Gamma_{ d \kappa } $,
$ \omega ( \pi \bullet , \bullet ) |_{\Gamma_{ d \kappa }}  $ is then 
given by $ \langle \Omega \bullet, \bullet \rangle $, where
\[ \Omega = 
  \left( \begin{array}{ll} \alpha & \beta - 1 
\\ \beta^t -1  & \gamma \end{array} \right)  \,, \]which proves the lemma.
\end{proof}

As is well known, and as will be seen 
in the proof of the next proposition, {\em any} locally defined 
Fourier Integral Operator can be represented by \eqref{eq:6.1}.
To compute its trace in terms of invariantly defined objects we also
have to recall the definition of the {\em Maslov index of a curve of
linear symplectomorphisms} -- see \cite{Da} for more details and
references. 

Thus let $ \Gamma ( t) \subset T^* \RR^n \times T^* \RR^n $, $
 a \leq t \leq b $,  be 
a curve of graphs of symplectomorphisms. Choose a subdivision
$ a = t_0  < t_1 < \cdots < t_k = b $, such that, for all 
$ j = 1 , \cdots , k $, there is a Lagrangian subspace $ M^j $
transversal to $ \Gamma( t ) $ and the diagonal, $ \Delta $, 
for $ t \in [ t_{ j -1} , t_j ] $. We now follow 
\cite{Da} and define the {\em Maslov index} of a
curve of linear symplectomorphisms as 
\begin{equation}
\label{eq:Maslov}
 \mu \stackrel{\rm{def}}{=} 
\frac12 \sum_{ j = 1}^k \left( s ( \Gamma ( t_{j-1} ), 
\Delta, M^j  ) - s ( \Gamma ( t_{ j } , \Delta, M^j ) ) \right) \,. \end{equation}
It is independent of the choice of the transversal Lagrangians, $ M^j $,
and of the subdivision.

We can now prove
\begin{prop}
\label{p:6.2}
Suppose that $ U (t)  $ is a family of 
Fourier Integral Operators 
defined using a family of pseudodifferential operators, $ A ( t) 
\in \Psi^{0,0} ( T^*  X) $,
as in \eqref{eq:2.2}:
\[  hD_t U( t) + A( t) U ( t)  = 0 \,, \ \ U(0) = U_0 \in \Psi_h^{0,0} ( 
T^* X ) \,. \]
Let us also assume that $ a_t =\sigma ( A ( t) ) $, the Weyl symbol 
(with a possible dependence on $ h$ in the subprincipal symbol part)
of $ A ( t) $, is real and generates a family of canonical transformations:
\begin{gather*}
 \frac{d}{dt} \kappa_t ( x, \xi ) =  ( \kappa_t )_* ( H_{a_t} 
( x, \xi )) \,, \ \ \kappa_0  ( x, \xi ) = ( x, \xi )\,,   \ \ 
( x, \xi ) \in T^* X  \\ 
\kappa_t ( 0 , 0 ) = ( 0 , 0) \,.
\end{gather*}
If 
\begin{equation}
\label{eq:non-d}
 \det ( 1 - d \kappa_T ( 0 , 0 ) ) \neq 0 \,, \end{equation}
then
\begin{gather}
\label{eq:6.3}
\begin{gathered}
\tr U (T) = i^{\mu ( T ) }
{ {(1+{\mathcal O}(h))e^{- i \int_0^T a_t ( 0 , 0 ) dt /h}} \over
{|{\det (d\kappa _T(0,0)-1)}|^{\frac12}}}   \sigma ( U_0 ) ( 0 , 0)   \,,
\end{gathered}
\end{gather}
where $ \mu ( T ) $ is the Maslov index of the curve of linear
symplectic transformations $ d \kappa_t ( 0 , 0 ) $, $ 0 \leq t \leq T $.
\end{prop}
\begin{proof}
Let us first assume that for $ 0 \leq t \leq T $ 
\begin{equation}
\label{eq:6.a}
  ( \kappa_t ( y, \eta ) ; (y , \eta)   ) \mapsto ( x (\kappa_t ( y , \eta)) 
 , x ) \ \text{is surjective near $ ( 0 , 0 )$} \,.
\end{equation}
We follow the presentation from \cite[Appendix a]{Gr}.
Let $a_t$ be the Weyl-symbol of $A_t$ defined modulo ${\mathcal
O}(h^2)$ (if there is a subprincipal symbol we include it
in the principal one and obtain an $h$ dependent symbol).
Consequently the influence of the subprincipal symbol
will be accounted for as an ${\mathcal O}(h)$-dependence in the
canonical transformation $ \kappa_t $.) Let $\kappa _t$ be the 
canonical transformation generated by $ H_{a_t } $ as described in the
statement of the proposition. We can then view $\kappa _t$ 
as the canonical transformation 
associated to $ U ( t) $ (defined modulo ${\mathcal O}(h^2)$) and 
we claim that 
\begin{equation}
\label{eqref:13}
U(t) u(x)={1\over (2\pi h)^n}\iint e^{{i\over h}(\phi
(t,x,\eta )-y\cdot \eta )}b(t,x,\eta ;h)u(y)dyd\eta \,, \end{equation}
where
\begin{equation}
\label{eq:14}
\partial _t\phi (t,x,\eta )+a_t(x,\partial _x\phi )=0,\ \ \ 
\phi (0,x,\eta )=x\cdot \eta \,, \end{equation}
so that $\kappa _t:\, (\partial _\eta \phi ,\eta
)\mapsto (x,\partial _x\phi )$. The amplitude $ b $  has to satisfy
\begin{gather}
\label{eq:15}
\begin{gathered}
(hD_t+ a_t^w(x,hD))(e^{i\phi (t,x,\eta )/h}b(t,x,\eta ;h)) =0\,, \\
(\partial _t\phi +hD_t+e^{-i\phi /h}a_t^w(x,hD)e^{i\phi
/h})(b)=0 \,. \end{gathered} \end{gather}
Here the Weyl symbol of $e^{-i\phi /h}a_t^w e^{i\phi /h}$
is $q_t(x,\xi )=a_t(x,\phi _x'+\xi )+{\mathcal O}(h^2)$, and
using that $\partial _t\phi =-a_t(x,\partial _x\phi )$,
we get
$$(hD_t+{\rm Op\,}(f_t(x,\xi )))b={\mathcal O}(h^2),$$
where $f_t(x,\xi )=a_t(x,\phi _x'+\xi )-a_t(x,\phi '_x)$
(and $\eta $ is just a parameter). This gives
\begin{equation}
\label{eq:16}
(hD_t+{1\over 2}\sum_1^n((\partial _{\xi
_j}a_t)hD_{x_j}+hD_{x_j}\circ (\partial _{\xi
_j}a_t)))b_0=0 \,, \end{equation}
for the leading part of $b$. With $\nu _t=\sum (\partial
_{\xi _j}a_t)\partial _{x_j}$, this can also be written
$$(\partial _t+\nu _t+{1\over 2}{\rm div\,}\nu
_t)b_0=0,$$ or
\begin{equation}
\label{eq:17}
(\partial _t+{\mathcal L}_{\nu _t})(b_0(t,x,\eta
)(dx_1\wedge ...\wedge dx_n)^{1/2})=0\,, \end{equation}
where  $ {\mathcal L}_{\nu _t} $ denotes the Lie derivative.
If we consider $b_0(t,x,\eta )(dx_1\wedge ...\wedge
dx_n)^{1/2}$ as a half-density on $\Lambda _{\phi
_{t,\eta }}=\{ (x,\partial _x\phi _t(x,\eta ))\}$, then
\eqref{eq:17} means that
\begin{equation}
\label{eq:18}
\kappa _t^*({{b_0(t,x,\eta )(dx_1\wedge ...\wedge
dx_n)^{1/2}}}|_{\Lambda _{\phi _{t,\eta
}}})={{(dx_1\wedge ...\wedge dx_n)^{1/2}}}|_
{\Lambda _{\phi _{0,\eta }}}\,.  \end{equation}

From \eqref{eq:5} it follows that the restriction of the
differential of $\kappa _t$ to $T\Lambda _{\phi _{0,\eta
}}$ followed by the $x$-space projection is given by
$$\delta _y\mapsto (\phi ''_{\eta x})^{-1}(\delta _y),$$
so \eqref{eq:18} says that
\begin{equation}
\label{eq:19}
{b_0(t,x,\eta )\over (\det \phi ''_{\eta x})^{1/2}}=1\,.
\end{equation}
We note that $ \det  \phi''_{\eta x} > 0 $ for $ 0 \leq t 
\leq T $.

From \eqref{eq:14} and $ d \phi_t ( 0,0 ) = 0$ (since 
$ \kappa_t ( 0, 0 ) = ( 0 , 0 ) $) we see that 
\[ \phi_T ( 0 , 0 ) = - \int_0^T a_t ( 0 , 0 ) dt \,. \]
Applying Lemma \ref{p:6.1} and Lemma \ref{l:6.1}
we obtain \eqref{eq:6.3}: under the assumption 
\eqref{eq:6.a} we only need one transversal Lagrangian in 
\eqref{eq:Maslov}, and we can take $ M $ from Lemma \ref{l:6.1}.
Then
\[ 
\begin{split}
\mu ( T ) & = 
\frac 12 \left( s ( \Gamma_{ d \kappa_0 } , \Delta , M ) - 
s ( \Gamma_{ d \kappa_T } , \Delta, M ) \right)  \\ 
& = \frac 12 \left( s ( \Delta , \Delta , M ) - 
s ( \Gamma_{ d \kappa_T } , \Delta, M ) \right) \\
& = - \frac12 s ( \Gamma_{ d \kappa_T } , \Delta, M ) =
\frac12  {\rm{sgn}} \; \left( \begin{array}{ll} \phi ''_{xx} &\phi
''_{x\eta }-1 \\ \phi ''_{\eta x}-1&\phi ''_{\eta \eta
}   \end{array} \right) \,. 
\end{split}
 \]
In the case \eqref{eq:6.a} does not hold for $ 0 \leq t \leq T $,
we have to choose different coordinates in which \eqref{eq:6.a} holds for 
$ t_{j-1} \leq t \leq t_j $. That gives  corresponding Lagrangians
$ M^j $ (defined as $ M $ was) and the phase shifts add up precisely to 
give \eqref{eq:Maslov}. In fact, we can conjugate $ U ( t ) $ by an
$h$-Fourier Integral Operator (so without affecting the trace), so that
for $ t_1 - \delta < t < t_2 + \delta $
the resulting operator is given by 
\[ \widetilde U ( t ) u ( x) 
={ i^\nu \over (2\pi h)^n}\iint e^{{i\over h}(\tilde \phi
(t,x,\eta )-y\cdot \eta )}\tilde b(t,x,\eta ;h)u(y)dyd\eta \,, \]
where we can arrange that $ \tilde b > 0 $, and that 
\eqref{eq:non-d} holds with $ T = t_1 , t_2 $.  We then use 
Lemma \ref{p:6.1} and the geometric discussion above to compute the trace:
\[ \tr  U ( t_1 ) = i^\nu i ^{ - \frac12 \tilde s } 
{ {(1+{\mathcal O}(h))e^{- i \int_0^{t_1} a_{t_1} ( 0 , 0 ) dt /h}} \over
{|{\det (d\kappa _{t_1} (0,0)-1)}|^{\frac12}}}   \sigma ( U_0 ) ( 0 , 0)  \,,
\]
where 
\[ \tilde s = - 
 {\rm{sgn}} \; \left( \begin{array}{ll} \tilde \phi ''_{xx} & \tilde \phi
''_{x\eta }-1 \\ \tilde \phi ''_{\eta x}-1& \tilde \phi ''_{\eta \eta
}   \end{array} \right) \,. \]
As in Lemma \ref{l:6.1} we interpret $ \tilde s $ as
\[ \tilde s = 
s ( \Gamma_{ d \kappa_{t_1} } , \Delta, M^2 )  \,, \]
where $ M_2 $ was chosen in new coordinates. Comparing it with the 
previous expression for the trace (where we put $ T = t_1 $ and $ 
M = M^1 $), we see that
\[ \nu = - s (  \Gamma_{ d \kappa_{t_1} } , \Delta, M^1 ) + 
s ( \Gamma_{ d \kappa_{t_1} } , \Delta, M^2 )  \,. \]
We now use $ \tilde U ( t ) $ to compute the trace at $ t = t_2 $
which in view of the expression for $ \nu $, and the fact that
$ s ( \Gamma_{ d \kappa_{0} } , \Delta, M ^1 ) = 0 $ is
\[ \begin{split}  \tr  U ( t_2) = & 
i^{ \frac12 (  s ( \Gamma_{ d \kappa_{0} } , \Delta, M ^1 ) 
- s (  \Gamma_{ d \kappa_{t_1} } , \Delta, M^1 ) + 
s ( \Gamma_{ d \kappa_{t_1} } , \Delta, M^2 )  - s 
( \Gamma_{ d \kappa_{t_2 }}, \Delta, M^2 ) )  }  \times \\
& 
{ {(1+{\mathcal O}(h))e^{- i \int_0^{t_1} a_{t_1} ( 0 , 0 ) dt /h}} \over
{|{\det (d\kappa _{t_1} (0,0)-1)}|^{\frac12}} }   \sigma ( U_0 ) ( 0 , 0)  \,,
\end{split}
\]
and comparing with \eqref{eq:Maslov} see that the power of $ i $ is given by 
the Maslov index for the curve $ d\kappa_{t} ( 0 , 0 ) $, $ 0 \leq t \leq
t_2 $. We can continue in the same way which gives us the final index
$ \mu ( T) $.

\end{proof}

We now want to evaluate the trace of $ M ( z, h ) ^k M' ( z , h ) $,
and for this we need to identify the Maslov factor and the phase. 
For this we recall the definition of the classical action:
\begin{equation}
\label{eq:action}
I( z) \stackrel{\rm{def}}{=} \int_{\gamma( z) } \xi d x \,. 
\end{equation}
The well known relation with the periods is given in 
\begin{lem}
\label{l:6.2}
Let $ q( z) $ and $q_{\circlearrowright} ( z ) $ be the
local time and the first return local time defined in \eqref{eq:3.3}
and \eqref{eq:3.3'}. Then 
\begin{gather*}
 (  q_{\circlearrowright}  ( z ) -q ( z )) |_{\gamma( z ) }
=  - \int_0^{T(z)} \sigma ( \partial_z P ( z) ) ( \exp
( t H_{p( z)} ) ( m_0(z)  ) dt \,, \\  
 (  q_{\circlearrowright}  ( z ) -q ( z )) |_{\gamma( z ) }
  = \frac{ d I }{dz } (z) \,. \end{gather*}
\end{lem}
\begin{proof}
The first identity follows directly from the definition and was 
already used in the proof of Lemma \ref{l:5.2}.

Since $ \partial_z p ( z) \neq 0 $, we can write $ p ( z) 
= c ( z) ( z - \tilde p ) $. Hence on $ p ( z) = 0 $, 
the equations for $ q $ and $ q_{\circlearrowright} 
$ are
\[  H_{\tilde p } q = - 1 \,, \ \ H_{\tilde p } q_{\circlearrowright} 
= -1 \,,  \]
and consequently $
 (  q_{\circlearrowright}  ( z ) -q ( z )) |_{\gamma( z ) }
  $ is the period of $ \gamma ( z) $, thought of as
an orbit of $ H_{\tilde p } $ on $ \tilde p = z $. We now introduce
an isotropic submanifold, $ \Gamma $,
 of $ T^* ( X \times \RR ) $, where the new 
variable (on $ \RR $) is denoted by $ \zeta $ with $ z $ its dual variable:
\[ \Gamma = \{  ( m ; ( \zeta , z ) ) \in T^* ( M \times \RR ) \; : \;
m \in \gamma ( z )\,, \ \zeta = q ( z) ( m ) \,,  z_0 \leq z \leq z_1 
\} \,.\]
The symplectic form $ d \xi \wedge dx + d z \wedge d \zeta $ 
vanishes on $ \Gamma $, and hence we obtain from Stokes's theorem:
\[\begin{split}
 I ( z_1 ) - I ( z_0) & = \int_{\gamma( z_1) } \xi dx - 
\int_{ \gamma ( z_2 ) } \xi d x  \\
& = \int_{ z_0}^{z_1} 
 \left( q ( z ) ( m_0 ( z ) ) - q_{\circlearrowright} ( z) 
( m_0 ( z ) ) \right) d z = \int_{z_0}^{z_1} T ( z) dz \,, 
\end{split}
\]
which proves the lemma.
\end{proof}

Using this lemma we will be able to identify the phase in the trace
of the monodromy operator. For that let
${\mathcal T}(z)$ be the
quantum time appearing in \eqref{eq:3.m}:
\[ {\mathcal T}(z) =  K ( z )^{-1} ( Q ( z ) - \QC ( z ) ) K ( z ) \,, \]
so that that formula becomes
\begin{equation}
\label{eq:6.22}
{
hD_zM(z)={\mathcal T}( z ) M(z) \,. 
}
\end{equation}
This and Proposition \ref{p:6.2} show that 
the phase factor in $ \tr M ( z , h )^k 
{\mathcal T}( z ) $ satisfies
\[ J_k ' ( z ) = k  
(  q_{\circlearrowright}  ( z ) -q ( z )) |_{\gamma( z ) }
 \,. \]
In fact, for any family of $ P ( z, h ) $ satisfying the assumptions of 
Proposition \ref{p:4.2}, we can associate to $ M ( z , h )^k $ (not
necessarily satisfying the non-degeneracy condition) 
a phase factor, $ J_k ( z ) $ which has to satisfy
\[ J_k ( z ) = k I ( z) + C_k \,, \]
We want to show that $ C_k =  0$. For that we note that if
we put $ P_\epsilon ( z , h ) = P ( z , h / \epsilon ) $ then 
the corresponding 
$ J _k( z) $ is 
 given by $ J_k ( z ) \epsilon $. 
On the other hand, the action corresponding to $ P _\epsilon $ is
$ k I ( z ) \epsilon $. Since we can consider $ P_\epsilon $ as 
another deformation of our operator we must have
\[ \ \forall \ \epsilon > 0 \,, \ \ 
  J_k ( z ) \epsilon = k I ( z ) \epsilon + C_k \,,\]
and for that we need that $ C_k = 0 $.


To obtain the Maslov factor we need to find
a family of symplectic transformations interpolating between the 
identity and the Poincar\'e map. For that let us fix $ z $ and
supress dependence on $ z $ in the subsequent formul{\ae}. We want 
to define a family $ M ( t ) : {\mathcal D}' ( \RR^{n-1} ) 
\rightarrow {\mathcal D}' ( \RR^{n-1} ) $, of $h$-Fourier Integral
Operators such that $ M ( 0 ) = Id $ and $ M ( T ) = M $.  To do
this we modify the definition of $ I_+ $ in  \eqref{eq:epm} to 
\[ I_+ (t) \; : \; {\rm{ker}}_{m_0 } P  
\ \longrightarrow  \  {\rm{ker}}_{ \exp t H_{  p } m_0 } ( P )  \]
We also  generalize the definition of $ K $ to 
\[ K ( t) \; : \; {\mathcal D}' ( \RR^{n-1} ) 
\ \longrightarrow \ 
{\rm{ker}}_{ \exp t H_{  p } m_0 } ( P ) \,, 
\]
defined using Proposition \ref{p:2.1} as in \eqref{eq:K}. 
We can now define 
\begin{equation}
\label{eq:mt}
M(t) \stackrel{\rm{def}}{=} K(t)^{-1} I_+ ( t ) K ( 0 ) 
\; : \; {\mathcal D}' (\RR^{n-1}) \ \longrightarrow {\mathcal D}' 
( \RR^{n-1} ) \,, 
\end{equation}
microlocally near $ ( 0, 0 ) $. This family has desired properties
and quantizes a curve of local symplectomorphism  $ \kappa_t $:
\[ \kappa_t = (\Phi_t)^{-1} \Psi_t
 \Phi_0 \; : 
\; T^* \RR^{n-1} \ \longrightarrow \ T^* \RR^{n-1} \,, \ \ 
\kappa_t ( 0, 0 ) = ( 0 , 0 ) \,, \]
where $ \Phi_t $ symplectically identifies a neighbourhood of $ ( 0 ,0 ) $ in 
$ T^* \RR^{n-1} $ with a submanifold of $ p = 0 $, 
$ S_t $, transversal to $ \gamma $
at $ \exp{ t H_{ p } } ( m_0 ) $, and 
$ \Psi_t : S_0 \longrightarrow S_t $ is the restriction of the flow 
$ \exp ( s H_{ p } ) $ to $ S_t $. The construction above allows
an arbitrary choice of $ S_t $ and  $ \Phi_t $.

We can summarize this discussion in 
\begin{prop}
\label{p:6.4}
Suppose that the orbit $ \gamma ( z) $ is primitive 
and $N$-fold non-degenerate in the sense that 
\eqref{eq:6.non} holds. Let $ I ( z) $ be the classical actions 
defined by \eqref{eq:action},  
and $ T (z )$ the periods of $ \gamma ( z)$.

If $ 
t $, $ 0 \leq t \leq T( z) $ parametrizes $ \gamma ( z ) $, let 
$ S_t $ be a family of submanifolds of $ p ( z)  = 0 $, transversal
to $ \gamma ( z) $ at $ t $, $ \Phi_t $ a symplectic identification of
$ S _t $ with a neighbourhood of $ ( 0 , 0 ) $ in $ T^* \RR^{n-1} $,
and $ \Psi_t : S_0 \rightarrow S_t $ the restriction of the 
flow to $ S_t $. Then for $ 0 < |k | \leq N $, 
\[  {\mathrm {tr}} M( z , h )^{k-1} hD_z M( z , h ) = 
\frac{ e^{ i k I( z) } e^{  i \nu_k  ( z) \frac{\pi}{2} }  
(  q_{\circlearrowright}  ( z ) -q ( z )) |_{\gamma( z ) }
}{
| (d C( z)_{m_0( z) })^k  - 1 |^{\frac12}} ( 1 + {\mathcal O} ( h ) ) \,, 
\]
where $ \nu_k ( z ) $ is the Maslov index of the curve of linear 
symplectic transformations:
\[  d ( \Phi_t^{-1} \Psi_t \Phi_0 )_{(0,0)} \,, \ \  0 \leq t \leq k T( z) 
\,. \]
\end{prop}

\medskip
\noindent
{\bf Remark.} 
The Maslov index $ \nu_k (z) $ is a locally constant function of 
$ z $: it does not change as long as \eqref{eq:6.non} holds. 
Its value may depend on the non-unique choices of the identifications
$ \Phi_t $, and the transversals $ S_t$. Since $ \exp ( i \nu_k  \pi/2 ) $
is determined uniquely (as it appears in the trace!), $ \nu $ is 
termined only modulo $ 4$. In the case when $ \gamma ( z ) 
\rightarrow \pi ( \gamma ( z )) $ is a diffeomorphism, 
with $ \pi : T^* X \rightarrow X $, the natural projection, 
a choice of transversals submanifolds in the base gives natural
$ S_t $'s in $ \{ p = 0 \} \subset T^* X $. Thus in the case of 
the geodesic flow on a Riemannian manifold $ \nu $ is the 
index of a closed geodesic.

The usual semi-classical trace formula for non-degenerate 
orbits follows from Theorem \ref{t:0} and the following
\begin{prop}
\label{p:6.3}
Suppose that the assumptions of Theorem \ref{t:0} are satisfied,
$ M ( z )$ is the quantum monodromy defined in \eqref{eq:3.mon},
and in addition the closed orbit $ \gamma = \gamma ( 0) $ 
is $ N$-fold non-degenerate
\eqref{eq:6.non}. Then, for $ k \ne 0 $, $ |k | \leq N $, and
$ g \in \CI_{\rm{c}}  ( \RR ) $, we have
\begin{equation}
\label{eq:p6}
\frac{1}{2 \pi i} 
{\rm{tr}} \; \int_{\RR} \widehat g  ( z / h ) M( z , h )^{k-1} \frac{d}{dz}
M( z , h ) \chi ( z ) d z 
= \frac{ e^{ i k S_\gamma/ h + i  \nu_{\gamma, k }  \frac{\pi}{2}   }
  T_\gamma
 g ( k T_\gamma ) }{ | \det ( (d C_\gamma)^k - I ) |^{\frac12} } + 
{\mathcal O} ( h ) \,,
\end{equation}
where $ T_\gamma $ is the primitive period of $ \gamma $, 
$ dC_\gamma $ is the
linear Poincar{\'e} map, $ S_\gamma $, the classical action of $ \gamma $, 
and $ \nu_{\gamma, k }  $ the Maslov index of $ k \gamma $.
\end{prop}
\begin{proof}
Since $ P$ is assumed to be 
self-adjoint, the subprincipal symbold of $ P $ is zero. Let
$\kappa _z$ be the Poincar{\'e} map and assume that $(d\kappa
_0(0,0))^k-1$ is non-degenerate. Let ${\mathcal T}(z)$ be the
quantum time appearing in \eqref{eq:6.22} above. 
Using the cyclicity of the trace, we can write the left hand side of
\eqref{eq:p6} as
\begin{equation}
\label{eq:6.23}
{
{1\over 2\pi }{\rm tr\,}\int_{\bf R}\chi
(z)\widehat{g}({z\over h})M(z)^k{\mathcal T}(z){dz\over h}. }
\end{equation}
The leading symbol of ${\mathcal T}(z)$ at the fixed point is
the period $T(z)=dI(z)/dz$, where $I(z)$ is the action
along the closed trajectory. Then to leading order, \eqref{eq:6.23}
becomes 
\begin{equation}
\label{eq:6.24}
{
{ i ^{ \mu } \over 2\pi }\int \chi (z)\hat{g }({z\over
h}){e^{ikI(z)/h}\over |\det((d\kappa
_z(0,0)^k-1)|^{\frac12} }T(z){dz\over h}, }
\end{equation}
where $ \mu $ is the Maslov index. Write
$E=z/h$, so that
$${I(z)\over h}={I(0)\over h}+I'(0)E+{\mathcal O}(h).$$
Then \eqref{eq:6.24} becomes, again to leading order,
\begin{equation}
\label{eq:6.25}
{{i ^ \mu \over 2\pi }\int \hat{g}(E)e^{ikI'(0)E}dE
{e^{ikI(0)/h}T(0)\chi (0)\over |\det((d\kappa
_0(0,0))^k-1)|^{\frac12} }}
{=
i^\mu {e^{ikI(0)/h}T(0)\chi (0) g (kT(0))\over
|\det((d\kappa
_0(0,0))^k-1)|^{\frac12} }.
}
\end{equation}

\end{proof}

The usual Gutzwiller trace formula for a more general class of operators
is given in 
\begin{thm}
\label{t:2}
Suppose that the assumptions of Theorem \ref{t:0} hold and that in 
addition $ \gamma $ is an $ N$-fold non-degenerate orbit in the 
sense that \eqref{eq:6.non} holds. Then in the notation of Proposition 
\ref{p:6.3} we have
\[
 {\mathrm{tr}} \; f ( P / h ) 
\chi ( P ) 
A =  \frac{1}{2\pi}   \sum_{k=-N}^N  
\frac{ e^{ i k S_\gamma/ h + i  \nu_{\gamma, k } \frac{\pi}{2}   }
  T_\gamma
\hat f  (-  k T_\gamma ) }{ | \det ( (d C_\gamma)^k - I ) |^{\frac12} } + 
{\mathcal O} ( h )    \,. \]
\end{thm}

%





\section*{Appendix}
\setcounter{section}{1}
\setcounter{equation}{0}
\setcounter{prop}{0}
\renewcommand{\thesection}{\Alph{section}}

In the classical treatment of  pseudo-differential operators, the
subprincipal symbol is invariant under coordinate changes when 
the pseudo-differential operators are considered as acting on 
half-densities. This invariance is particularly nice in the
Weyl calculus, where the subprincipal symbol is contained in the 
leading symbol -- see \cite[Sect.18.5]{Hor}.

For the reader's convenience we present here a self-contained discussion
of the analogous result in the semiclassical setting.

We use the informal notation for sections of the half-density bundles:
\begin{gather*}
u \in \CI ( X , \Omega_X^{\frac12}) \ \ \Longleftrightarrow \ \ 
u = u ( x ) |dx|^{\frac12}  \,, \\
a \in S^{ 0, 0} ( T^*X , \Omega_{T^*X }^{\frac12} ) \ \ 
\Longleftrightarrow \ \ 
a = a ( x, \xi ) |dx|^{\frac12} |d \xi |^{\frac12}  \,, \end{gather*}
which captures the transformation laws under changes of coordinates:
\[ u ( x ) |d x |^{\frac12} =  \tilde u ( \tilde x ) |d \tilde x|^{\frac12}\,,
\  \tilde x = \kappa ( x ) \ \Longleftrightarrow \ 
\tilde u ( \kappa ( x ) ) | \kappa ' ( x ) |^{\frac12 } = u ( x  ) \,,
\]
where for a linear tranformation $ A $ we denote its determinant by 
$ | A | $.

We observe that the half-density sections over $ T^* X $ are identified with 
functions if we consider symplectic changes of variables, and in particular
\begin{equation}
\label{eq:kappa}
 ( x , \xi ) \longmapsto ( \tilde x, \tilde \xi ) = 
( \kappa ( x ) \,, ( \kappa' ( x ) )^{t} \xi ) \,.
\end{equation} 

As stated after \eqref{eq:weyl} in this paper we considered pseudodifferential
operators acting on half-densities:
\[ \Psi_h ^{ m , k } ( X ) = \Psi_h ^{m, k } ( X , \Omega_X^{\frac12} ) \,,\]
with distributional kernels given by 
\begin{equation}
\label{eq:weyla}
  K_a ( x , y ) | dx|^{\frac12} | d y |^{\frac12} = 
 \frac{1}{ ( 2 \pi h )^n } 
\int   a \left( \frac{x + y }{2}  , \xi \right
) e^{ i \langle x -  y, \xi \rangle / h }  d \xi |dx|^{\frac12}
|dy |^{\frac12}   \,. \end{equation}

We will show that 
\begin{equation}
\label{eq:inva}
  K_a ( x , y ) | dx|^{\frac12} | d y |^{\frac12} =
  K_{\tilde a} ( \tilde x , \tilde y ) | d \tilde x|^{\frac12} 
| d \tilde y |^{\frac12} \ \ \ \text{with  $ \tilde a ( \tilde x , \tilde 
\xi ) = 
a ( x , \xi ) + {\mathcal O} ( h^2 \langle \xi \rangle^{-2} ) $,}
\end{equation}
where  $ ( \tilde x , \tilde \xi ) $ is given by \eqref{eq:kappa}.

To establish \eqref{eq:inva} 
we start with its right hand side using 
the coordinates
$ \tilde x = \kappa ( x) $ and $ \tilde y = \kappa ( y ) $:
\[ \frac{1}{ ( 2 \pi h )^n } 
\int   \tilde a \left( \frac{\tilde x + \tilde y }{2}  , \tilde \xi \right
) e^{ i \langle \tilde x -  \tilde y, \tilde 
\xi \rangle / h }  d \tilde \xi |d \tilde x|^{\frac12}
|d \tilde y |^{\frac12}   \,. \]
Making a substition we obtain
\[ \frac{1}{ ( 2 \pi h )^n } 
\int   \tilde a \left( \frac{ \kappa (  x ) + \kappa(  y)  }{2}  , \tilde 
\xi \right
) e^{ i \langle \kappa(  x)  -  \kappa(  y) , \tilde 
\xi \rangle / h }  d \tilde \xi  | \kappa ' ( x ) | ^{\frac12}
| \kappa ' ( y ) |^{\frac12} 
|d  x|^{\frac12}
|d  y |^{\frac12}    \,. \]
We now apply the ``Kuranishi trick'' and for that write 
\begin{gather}
\label{eq:kur}
\begin{gathered}
 \kappa ( x) - \kappa ( y ) =  F ( x , y ) ( x - y ) \,, \ \ 
F ( x, y ) = \kappa' \left( \frac{x+y}{2} \right) + 
{\mathcal O} ( ( x-y)^2 ) \,, \\ 
\kappa( x  ) + \kappa ( y ) = \kappa \left( \frac{x+y}{2} \right) + 
{\mathcal O} ( ( x - y ) ^2 ) \,.
\end{gathered}
\end{gather}
We put $ \xi = F ( x , y )^{t} \tilde \xi $ and
rewrite the expression above as
\[
\begin{split}
&  \frac{1}{ ( 2 \pi h )^n } 
\int  \left(  \tilde a \left( \kappa \left(
\frac{   x +   y  }{2} \right)  ,  ( K( x, y )^t)^{-1} \xi \right) 
+ {\mathcal O} ( x - y )^2 \right) 
 e^{ i \langle   x  -    y ,  \xi 
\rangle / h }  | K( x , y )^t  |^{-1} \\
& \ \ \ \ \ \ \ \ \ \ \ \ \ d  \xi  | \kappa ' ( x ) | ^{\frac12}
| \kappa ' ( y ) |^{\frac12} 
|d  x|^{\frac12}
|d  y |^{\frac12}   = \\
&  \frac{1}{ ( 2 \pi h )^n } 
\int  \left(  \tilde a \left( \kappa \left(
(   x +   y  )/{2} \right) ,   \left(  \kappa' \left( 
(x+ y ) /2 \right)^t \right) ^{-1} \xi \right) 
+ {\mathcal O} ( x - y )^2 \right) 
 e^{ i \langle   x  -    y ,  \xi 
\rangle / h } | \kappa'(  ( x + y )/2 ) |^ {-1} \\
& \ \ \ \ \ \ \ \ \ \ \ \ \ d  \xi  | \kappa ' ( x ) | ^{\frac12}
| \kappa ' ( y ) |^{\frac12} 
|d  x|^{\frac12}
|d  y |^{\frac12}  \,,  \end{split} \]
and the terms $ {\mathcal O} ( (x-y)^2 ) $ will  influence the 
symbol only modulo $ {\mathcal O} (\langle \xi \rangle^{-2} h^2) $ 
(by integration by parts based on $ (x - y ) \exp ( \langle x - y , \xi 
\rangle / h ) = h D_\xi \exp ( \langle x - y , \xi 
\rangle / h ) $ ), and hence can be neglected.

We now observe that
\[  | \kappa'(  ( x + y )/2 ) |^2  = | \kappa ' ( y ) |  | \kappa ' ( x ) |
+ {\mathcal O} ( ( x - y ) ^2 ) \,, \]
and  consequently 
\[ \begin{split} & K_{\tilde a } ( \tilde x , \tilde y ) = \\
& \ \  \frac{1}{ ( 2 \pi h )^n } 
\int  \left(  \tilde a \left( \kappa \left(
(   x +   y  )/{2} \right) ,   \left(  \kappa' \left( 
(x+ y ) /2 \right)^t \right) ^{-1} \xi \right) 
+ {\mathcal O} ( h^2 \langle \xi \rangle^{-2} )  \right) 
 e^{ i \langle   x  -    y ,  \xi 
\rangle / h } | d  \xi  |d  x|^{\frac12}
|d  y |^{\frac12}  \,, \end{split}
\]
which is the same as $ K_a ( x , y ) |d x |^{\frac12} | d y |^{\frac12} 
$, if
\[  a ( x , \xi ) = \tilde a \left( \kappa \left(
x \right) ,   \left(  \kappa' \left( 
x  \right)^t \right) ^{-1} \xi \right) + {\mathcal O} ( h^2 \langle 
\xi \rangle ^{-2 } ) \,. \]
This proves \eqref{eq:inva} completing the appendix.

\end{document}